\documentclass[graybox]{svmult}
\usepackage{type1cm}        % activate if the above 3 fonts are
\usepackage{makeidx}         % allows index generation
\usepackage{graphicx}        % standard LaTeX graphics tool
\usepackage{subfig}          % for subfigures...               
\usepackage{multicol}        % used for the two-column index
\usepackage[bottom]{footmisc}% places footnotes at page bottom

\usepackage{newtxtext}       % 
\usepackage[varvw]{newtxmath}% selects Times Roman as basic font

\usepackage{algpseudocode} 
\usepackage{algorithm} 
\usepackage{amsmath}
\usepackage{diagbox}
\newcommand*{\C}{\mathbb{C}}%.............................C
\newcommand*{\R}{\mathbb{R}}%.............................R
\newcommand{\vect}[1]{\mbf{#1}} % vectorial notation
\renewcommand{\d}{\ \mathrm{d}} % d symbol inside integrals
\newcommand{\dt}{\ \mathrm{dt}} % dt (time increment)
\DeclareMathOperator{\spann}{span}

\makeindex             % used for the subject index
\begin{document}

\title{Factorizations and fast diagonalization for the heat equation}
% Use \titlerunning{Short Title} for an abbreviated version of
% your contribution title if the original one is too long
\author{A. Bressan\orcidID{0000-0003-3730-668X} and\\ A. Kushova\orcidID{0000-0002-5243-2249} 
and\\ G. Loli\orcidID{0000-0002-2192-3889} and\\ M. Montardini\orcidID{0000-0003-0729-063X} and\\ G. Sangalli\orcidID{0000-0002-5642-1969} and\\ M. Tani\orcidID{0000-0002-7070-6776}}
% Use \authorrunning{Short Title} for an abbreviated version of
% your contribution title if the original one is too long
\institute{Andrea Bressan \at IMATI-CNR, Via Ferrata 5/a, 27100 Pavia, Italy,\\ \email{andrea.bressan@imati.cnr.it}
\and Alen Kushova \at  Dipartimento di Matematica, University of Pavia, Via Ferrata 5, 27100 Pavia, Italy \\ \email{alen.kushova@unipv.it}
\and Gabriele Loli \at Dipartimento di Matematica, University of Pavia, Via Ferrata 5, 27100 Pavia, Italy \\ \email{gabriele.loli@unipv.it}
\and Monica Montardini \at Dipartimento di Matematica, University of Pavia, Via Ferrata 5, 27100 Pavia, Italy \\ \email{monica.montardini@unipv.it}
\and Giancarlo Sangalli \at Dipartimento di Matematica, University of Pavia, Via Ferrata 5, 27100 Pavia, Italy \\ \email{giancarlo.sangalli@unipv.it}
\and Mattia Tani \at Dipartimento di Matematica, University of Pavia, Via Ferrata 5, 27100 Pavia, Italy \\ \email{mattia.tani@unipv.it}
}
%
% Use the package "url.sty" to avoid
% problems with special characters
% used in your e-mail or web address
%
\maketitle

\abstract*{This work investigates diagonalization-based methods for efficiently solving linear evolution problems, with a particular focus on the heat equation. 
The plain diagonalization of the differential operator, though effective for elliptic problems where fast diagonalization can be used, exhibits instability when applied to the heat equation. 
To address this difficulty, we examine three alternative approaches, based on LU~factorization, a suitable arrowhead factorization, and a low-rank modification. These methods introduce more robust factorizations of the time derivative, ensuring both computational efficiency and stability. 
}

\abstract{This work investigates diagonalization-based methods for efficiently solving linear evolution problems, with a particular focus on the heat equation. 
The plain diagonalization of the differential operator, though effective for elliptic problems where fast diagonalization can be used, exhibits instability when applied to the heat equation. 
To address this difficulty, we examine three alternative approaches, based on LU~factorization, a suitable arrowhead factorization, and a low-rank modification. These methods introduce more robust factorizations of the time derivative, ensuring both computational efficiency and stability. 
}

\def\gd{\ensuremath{\Omega_*}}%{\ensuremath{\mathcal{Q}}}
\def\std{\ensuremath{\Omega}}%{\ensuremath{\mathcal{Q}}}
\def\sd{\ensuremath{\Omega_s}}%{\ensuremath{\Omega}}
\def\td{\ensuremath{\Omega_t}}%{\ensuremath{\mathcal{I}}}
\def\to{\rightarrow}
\def\mbf#1{\mathbf{#1}}
\def\over#1{\ensuremath{\overset{\text{\tiny $\leftrightarrow$}}{\mbf{#1}}}}
\def\band{\ensuremath{\mathscr{b}}} % band width
\def\rank{\ensuremath{\mathscr{r}}} % rank of low rank modification

\section{Introduction}\label{sec:intro} 
Solutions of linear evolution problems are typically computed numerically by sequential time stepping methods.
Maday and Rønquist \cite{maday2008parallelization} proposed a method to parallelize this sequential computation by diagonalizing the time derivative matrix and then solving all time steps concurrently resulting in a Parallel-in-Time (PinT) direct solver.
Their strategy has since been further investigated and an overview of several variants, collectively known as \emph{ParaDiag} algorithms, is provided in \cite{gander2020paradiag}.  
This approach relies on the Cartesian structure of the space-time domain and the tensor product structure of the discretization space, which induces a Kronecker structure in the differential operator.
In this it is similar to earlier methods like \emph{alternating direction implicit} (ADI) \cite{birkhoff1962alternating} and \emph{fast diagonalization} (FD) \cite{tensormethods1964}, originally developed for efficiently solving finite difference discretizations of the Poisson problem.
This chapter recalls and compares three methods exploiting the Kronecker structure for the heat equation and highlighting their common structure that is also shared by the original approach of Maday and Rønquist.
These methods primarily differ in how they handle the inversion of the time derivative operator. 

The model problem considered is the heat equation. % and its space-time discretization.
Assuming homogeneous Dirichlet boundary and initial conditions, it is formulated as follows: find $u$ such that 
\begin{equation}
\label{eq:model_problem}
%\text{find } u \text{ such that }
    \begin{cases}
        \partial_t u - \Delta_s u  = f& \quad  \text{in } \,\phantom{\partial}\sd \times \td,\\
        u   = 0& \quad  \text{on } \partial\sd \times \td,\\
        u   = 0& \quad  \text{in }  \,\phantom{ \partial}\sd \times \{ 0 \},
    \end{cases}
\end{equation} 
where $\sd \subset \R^d$ is a bounded Lipschitz domain, with boundary $\partial\sd$, and $\td = [0,T]\subset \R$ is the time interval.   For $f\in L^2([0,T]\times \Omega_s)$  problem~\eqref{eq:model_problem} has a unique weak solution $u$ in the Bochner space of functions in $L^2(0,T; H^1_0(\Omega_s))$ with time derivative in $L^2(0,T;H^{-1}(\Omega_s))$ satisfying the initial condition, see \cite{evans-book} for details.
% Here, $\gamma>0$ is the heat capacity constant, $\nu > 0$ is the thermal conductivity constant and $f$ denotes the source term.
% For the sake of simplicity we set $\gamma = \nu = 1$. 

 System~\eqref{eq:model_problem} is a linear problem of the form $A u=f$ where the operator $A$ exhibits a tensor product structure:
\begin{equation}\label{eq:problem-continuous}
      A=  A_t \otimes M_s + M_t \otimes A_s,
\end{equation}
with $M_s$ and $M_t$ the identity operators seen as the Riesz operators of the $L^2$-scalar product. 
Here, $A_t$ is the time derivative, and $A_s$ is the space-Laplacian operator, which is elliptic and $M_s$-self-adjoint.
Consequently, a generalized eigendecomposition of the pair $(A_s,M_s)$ is possible and stable, which forms a crucial assumption in the methods discussed in this chapter.

If such a problem is discretized with a finite dimensional space $V=V_s\otimes V_t$, the structure is preserved leading to the linear system
\begin{equation}\label{eq:kronecker-structure-discrete}
    \mbf A \mbf u = (\mbf A_t \otimes \mbf M_s+\mbf M_t \otimes \mbf A_s)\mbf u= \mbf f,
\end{equation}
where the correspondence with the continuous operators is given by matching letter. 
Introducing $N_\square=\dim V_\square$ and $N=\dim V=N_sN_t$, then
$\mbf u, \mbf f\in \R^N$, $\mbf A\in\R^{N\times N}$, $\mbf M_s, \mbf A_s \in \R^{N_s \times N_s}$, and $\mbf M_t, \mbf A_t \in \R^{N_t \times N_t}$.

Efficient solvers for problems of the form \eqref{eq:kronecker-structure-discrete},
are well established in literature (see  \cite{simoncini2016computational} for a survey of methods requiring loose assumptions on the Kronecker factors).
In this chapter, we consider methods that are based on diagonalization and to illustrate the core concept we assume that
the generalized eigendecompositions of both pairs $(\mbf{A}_t, \mbf{M}_t)$ and $(\mbf{A}_s, \mbf{M}_s)$ are available: 
this means that there are diagonal matrices $\mbf\Lambda_t \in \C^{N_t\times N_t}$, $\mbf\Lambda_s \in \C^{N_s\times N_s}$ and invertible matrices $\mbf U_t \in \C^{N_t\times N_t},\mbf U_s \in \C^{N_s\times N_s}$ such that
\begin{align}
\label{eq:eigen-decomposition-t}
\mbf{A}_t \mbf{U}_t & = \mbf{M}_t \mbf{U}_t \mbf{\Lambda}_t,
\\\label{eq:eigen-decomposition-s}
\mbf{A}_s \mbf{U}_s & = \mbf{M}_s \mbf{U}_s \mbf{\Lambda}_s. 
\end{align}
By introducing the matrices $ \widetilde{\mbf{U}}_t = \left( \mbf{M}_t \mbf{U}_t \right)^{-1}$ and $ \widetilde{\mbf{U}}_s = \left( \mbf{M}_s \mbf{U}_s \right)^{-1}$, we obtain the factorizations
\begin{align*}
\mbf{A}_t &= \widetilde{\mbf{U}}_t^{-1} \mbf{\Lambda}_t \mbf{U}_t^{-1}, && \mbf{M}_t  = \widetilde{\mbf{U}}_t^{-1} \mbf{I}_{N_t} \mbf{U}_t^{-1},
\\\mbf{A}_s &= \widetilde{\mbf{U}}_s^{-1} \mbf{\Lambda}_s \mbf{U}_s^{-1}, && \mbf{M}_s  = \widetilde{\mbf{U}}_s^{-1} \mbf{I}_{N_s} \mbf{U}_s^{-1},
\end{align*}
with $\mbf{I}_{\square} \in \R^{\square\times \square}$ %and $\mbf{I}_{t} \in \R^{N_t \times N_t}$ 
the identity matrices. 
%In the methods presented, it will not be necessary to invert $\mbf M_t\mbf U_t$ and $\mbf M_s\mbf U_s$, as the definitions of $\widetilde{\mbf U}_s$ and $\widetilde{\mbf U}_t$ will be provided in closed form.
The above expressions lead to the following factorization of $\mbf A$ 
\begin{equation}\label{eq:A-factorizaion-simplified}
\mbf{A} = (\widetilde{\mbf{U}}_t\otimes \widetilde{\mbf{U}}_s)^{-1}(\mbf{\Lambda}_t\otimes \mbf{I}_{N_s} + \mbf{I}_{N_t}\otimes \mbf{\Lambda}_s)(\mbf{U}_t\otimes\mbf{U}_s)^{-1}.
\end{equation}
The factorization \eqref{eq:A-factorizaion-simplified} permits to solve \eqref{eq:kronecker-structure-discrete} by applying the inverse of each factor in turn. This approach is here called \emph{Diagonalization in Time} (DT) method. 
Furthermore, the Kronecker structure of $\widetilde{\mbf{U}}_t\otimes \widetilde{\mbf{U}}_s$ and ${\mbf{U}}_t\otimes {\mbf{U}}_s$ allows for matrix-vector multiplication with a reduced computational cost of 
$O(N_sN + N_t N)$, in contrast to the $O(N^2)$ complexity, and $(\mbf{\Lambda}_t\otimes \mbf{I}_{N_s} + \mbf{I}_{N_t}\otimes \mbf{\Lambda}_s)$ is diagonal allowing for inversion at a cost of $O(N)$ operations.
This highlights how the tensor/Kronecker structure facilitates the solution of the system. However, it hides important details about the
existence, stability and computational cost of the generalized eigendecompositions~\eqref{eq:eigen-decomposition-t}  and \eqref{eq:eigen-decomposition-s}.

The existence and stability of the eigendecomposition of the pair $(\mbf A_s, \mbf M_s)$, i.e. \eqref{eq:eigen-decomposition-s}, 
%is assumed for our purposes and it can be proved for the case of the heat equation.
follows for the heat equation from $\mbf A_s$ and $\mbf M_s$ being symmetric and positive definite. In such a case 
%Indeed we can impose the generalized eigenvector matrix is orthonormal with respect to $\mbf{M}_s$.  
\begin{equation} \label{eq:orthogonality} 
    \widetilde{\mbf{U}}_{s} = \left( \mbf{M}_s \mbf{U}_s \right)^{-1} =  \mbf{U}_{s}^T \quad \text{i.e.} \quad \mbf{U}^T_{s} \mbf{M}_{s} \mbf{U}_{s} = \mbf{I}_{N_s}. 
\end{equation}
Consequently, the condition number of the eigenvector matrix $\mbf{U}_s$ satisfies
$$ \kappa_2 \left( \mbf{U}_{s} \right) :=  \|\mbf{U}_{s}\|_2 \|\mbf{U}_{s}^{-1}\|_2 = \sqrt{\kappa_2(\mbf{M}_{s})},$$
which is often small or moderate. %, that is a requirement for the stability of the method. 

%However, the main challenge remains the computational cost which is typically $O(N_s^3)$ and results prohibitively high.
% However, this procedure relies on the computation of factorizations \eqref{eq:eigen-decomposition-t} and \eqref{eq:eigen-decomposition-s}. 
Assuming $N_s > N_t$, which is typically the case when $d > 1$, the computational cost of this step is in general $O(N_s^3)$ (see e.g. \cite[Chapter 8]{Golub2012}), and results prohibitively high.
% This limits the applicability of the described methods to two specific scenarios.
% In the first scenario, 
However, the generalized eigendecomposition of the pair $(\mbf{A}_s,\mbf{M}_s)$ can be computed with a reduced cost % $O( N_s^{\epsilon})$, where $\epsilon<3$. 
% This occurs, for example, 
when $\Omega_s$ is a Cartesian domain. In such a case Fast Diagonalization (FD) can be performed, see Section \ref{sec:fast-diagonalization}. 
%In this case the methods act as fast direct solvers. 

Another interesting case is when   a generalized eigendecomposition with reduced cost is possible for a suitable approximation 
\begin{equation}\label{eq:modified-matrix-preconditioner}
\widehat{\mbf A} = \mbf{A}_t \otimes \widehat{\mbf{M}}_s + \mbf{M}_t \otimes \widehat{\mbf A}_s,
\end{equation}
of $\mbf A$, that is one can use FD on the pair $(\widehat{\mbf{A}}_s,\widehat{\mbf{M}}_s)$. In such a case $\widehat{\mbf{A}}$ has the role of a preconditioner.
Both cases are exemplified in the literature and have been applied in different
contexts, including finite differences \cite{tensormethods1964, Buzbee1970}, finite elements \cite{Patera1986, Bjontegaard2009}, spectral elements \cite{Couzy1994, Lottes2005, Bjontegaard2009} and isogeometric analysis \cite{Montardini2018}.
Section~\ref{sec:heat_equation} reports numerical result of both cases
applied to spline based Galerkin discretizations.
%Numerical results for both cases are presented in Section~\ref{sec:heat_equation}.
%These approaches have been applied in different contexts, including finite differences \cite{tensormethods1964, Buzbee1970}, finite elements \cite{Patera1986, Bjontegaard2009}, spectral elements \cite{Couzy1994, Lottes2005, Bjontegaard2009} and isogeometric analysis \cite{Montardini2018}.

Concerning \eqref{eq:eigen-decomposition-t}, observe that the time derivative, at the continuous level, does not admit eigenfunctions fulfilling the homogeneous initial condition. Therefore, at the discrete level, the generalized eigendecomposition~\eqref{eq:eigen-decomposition-t} of the pair $(\mbf A_t,\mbf M_t)$, when computable, 
should be seen as a numerical artifact and indeed becomes instable in some cases, e.g. for high-degree spline based Galerkin discretization.
These instabilities are exemplified for both finite differences and Galerkin methods in Section~\ref{sec:diagonalization-in-time}, and they
may affect the accuracy of the solution.

This chapter presents three alternative methods, here called: \emph{LU~factorization} (LU)
\cite{Haidvogel1979, Shen1994, Wang2013, bressan2023space}, \emph{arrowhead factorization} (AR), originally  proposed in the context of isogeometric analysis \cite{MR2152382}, %as a preconditioner on parametrized domains that is a direct solver on the parametric domain in \cite{loli2020efficient} 
and \emph{low-rank modification} (LR), here introduced similarly to the methods in \cite{kressner2023improved,gander2024new}. %\marginpar{[6,5] sono simili ma non uguali perchè non diagonalizzano in spazio.}
Among them, AR and LR are space-time parallelizable, while LU is only parallelizable in space. 
Each method corresponds to a factorization of $\mbf A$, that replaces \eqref{eq:A-factorizaion-simplified}, of the form:
\begin{equation}\label{eq:generic-factorization}
    \mbf{A} = (\widetilde{\mbf{U}}_t\otimes \widetilde{\mbf{U}}_s)^{-1} \mbf{T}(\mbf{U}_t\otimes\mbf{U}_s)^{-1},
\end{equation}
where $\mbf U_s$ and $\widetilde{\mbf{U}}_s$ fulfill \eqref{eq:eigen-decomposition-s}, while $\mbf U_t$, $\widetilde{\mbf{U}}_t$ and $\mbf{T}$ are method specific.
Then, \eqref{eq:kronecker-structure-discrete} is solved following the template described by Algorithm \ref{alg:gen_fast_diagonalization}.

%Then, $\mbf A$ can be inverted as in Algorithm~\ref{alg:gen_fast_diagonalization} that is the template of the presented methods.
\begin{algorithm}[H]
    \caption{Common structure for the methods DT, LU, AR, LR.}\label{alg:gen_fast_diagonalization}
    \begin{algorithmic}[1]
    \State compute $\mbf{U}_s,\widetilde{\mbf{U}}_s,\mbf{\Lambda}_s$ as in \eqref{eq:eigen-decomposition-s}% \hfill with cost $O(\sum_{l=1}^d N_{s,l}^3)$
    \State compute $\mbf{U}_t,\widetilde{\mbf{U}}_t,\mbf{T}$
    \State compute $\mbf{y} = (\widetilde{\mbf{U}}_t\otimes\widetilde{\mbf{U}}_s) \mbf{f}$ \hfill
    \State solve $ \mbf{T}\mbf{z} = \mbf{y}$
    \State compute $\mbf{u} = (\mbf{U}_t\otimes\mbf{U}_s) \mbf z$ 
    \end{algorithmic}
\end{algorithm}
Notice that, when the methods are used to apply a preconditioner, 
steps 1--2 of Algorithm \ref{alg:gen_fast_diagonalization} are computed only once as the setup phase, while steps 3–5 are repeated for each application of the preconditioner and form the application phase.

In what follows, Section~\ref{sec:fast-diagonalization} describes the
FD method. Section~\ref{sec:diagonalization-in-time} describes the diagonalization in time method  introduced in \cite{maday2008parallelization} and reports specific examples of numerical instabilities 
for both finite difference and Galerkin discretizations of the heat equation.
The stable methods  LU, AR and LR are described in details in Section~\ref{sec:stable_time_factorizations} and numerical results are reported in Section~\ref{sec:heat_equation}.
Finally in Section \ref{sec:conclusions} the theoretical and numerical results are summarized. %\marginpar{Da rivedere in base alle modifiche sezione 1.2}

\section{Fast diagonalization}\label{sec:fast-diagonalization}
The generalized eigendecomposition in space \eqref{eq:eigen-decomposition-s} has computational cost of order $O(N_s^3)$, 
that is too demanding to be practical. However, it can be achieved at a reduced cost when the spatial domain has an inner Cartesian structure that is reflected in the function space
\begin{equation}\label{eq:space-tensor-product}
\Omega_s=\Omega_{s,1}\times\dots\times\Omega_{s,d},\qquad  V_s=V_{s,1}\otimes\dots\otimes V_{s,d},\end{equation}
in which case
\begin{align}
\label{eq:Ms-tp}\mbf M_s &= \mbf M_{s,d} \otimes \ldots \otimes \mbf M_{s,1},
\\\label{eq:As-tp}\mbf A_s &= \sum_{j=1}^d \mbf{K}_{j,d} \otimes \ldots \otimes \mbf{K}_{j,1},\quad\text{where}\quad  \mbf{K}_{j,l}=\begin{cases} \mbf{A}_{s,l} & \text{for }l=j,\\\mbf{M}_{s,l}& \text{for }l\ne j,\end{cases}
\end{align}
with $\mbf M_{s,l},\mbf A_{s,l} \in \R^{N_{s,l} \times N_{s,l}}$ for $l=1,\ldots,d$ representing respectively the identity and the Laplace operator on $V_{s,l}$.
The FD method applies more generally when $\mbf A_s$ is a sums of
kronecker matrices of the form 
    \begin{equation}\label{eq:kronecker_structure_of_stiffness_matrix}
        \mbf A_s = \sum_{j=1}^m c_j \mbf{K}_{j,d} \otimes \ldots \otimes \mbf{K}_{j,1},
    \end{equation} 
    with $c_j\in\R$, $\mbf K_{j,l} \in \{\mbf M_{s,l},\mbf A_{s,l}\}\subset\R^{N_{s,l} \times N_{s,l}}$ such that  $(\mbf A_{s,l}, \mbf M_{s,l})$ admits a generalized eigendecomposition, for all $l=1,\dots,d$.

In this case, for $l=1,\ldots,d$ there are $\mbf{U}_{s,l}$, $\mbf{\Lambda}_{s,l}$
that satisfy
\[
\mbf{A}_{s,l} \mbf{U}_{s,l} = \mbf{M}_{s,l} \mbf{U}_{s,l} \mbf{\Lambda}_{s,l}.
\]
Then, the generalized eigendecomposition of 
the pair $(\mbf{A}_s,\mbf{M}_s)$ of equation~\eqref{eq:eigen-decomposition-s}
 is realized with
\[
\begin{aligned}
    \mbf{U}_{s} &:=\mbf{U}_{s,d}\otimes\dots\otimes\mbf{U}_{s,1},\\
    \mbf \Lambda_s &:= \sum_{j=1}^m c_j\mbf{L}_{j,d} \otimes \ldots \otimes \mbf{L}_{j,1},\quad\text{where}\quad \mbf{L}_{j,l}=\begin{cases}
\mbf{I}_{N_{s,l}} & \text{for }\mbf{K}_{j,l}=\mbf{M}_{s,l},
\\ \mbf{\Lambda}_{s,l} & \text{for }\mbf{K}_{j,l}=\mbf{A}_{s,l}.\end{cases}
\end{aligned}
\]
Computing the generalized eigendecomposition of 
the pair $(\mbf{A}_{s,l},\mbf{M}_{s,l})$ has a computational cost of order 
$O(N_{s,l}^3)$, so that Step~1 in Algorithm~\ref{alg:gen_fast_diagonalization}
has a cost of $O(\sum_{l=1}^dN_{s,l}^3)$.
Notice that, assuming $N_{s,l}=N_s^{1/d}$, for $l=1,\dots,d$, then the eigenvalues are computed with cost $O(N_{s}^{3/d})$, that for $d=2$ is $O(N_s^{3/2})$, and for $d=3$ is $O(N_s)$.

 When the domain is not Cartesian, but it is parametrized over a Cartesian domain one can construct a preconditioner based on the technique above as exemplified in Section~\ref{sec:numeric-preconditioner}.

\section{Diagonalization in time  and numerical instabilities}\label{sec:diagonalization-in-time}
This method was proposed in \cite{maday2008parallelization} and is based on the generalized eigendecomposition \eqref{eq:eigen-decomposition-t}, further combined  with the generalized eigendecomposition \eqref{eq:eigen-decomposition-s}.
%This method is obtained when the generalized eigendecomposition \eqref{eq:eigen-decomposition-t} of the time derivative can be computed \cite{maday2008parallelization} and combined with the generalized eigendecomposition \eqref{eq:eigen-decomposition-s}. 
In this case, the factorization \eqref{eq:generic-factorization} reduces to %a plain diagonalization in space-time of the differential operator
\eqref{eq:A-factorizaion-simplified}, i.e
\[
\mbf{T} = \mbf{\Lambda}_t\otimes \mbf{I}_{N_s} + \mbf{I}_{N_t}\otimes \mbf{\Lambda}_s,
\] 
and the inverse of $\mathbf{A}$ can be applied with Algorithm~\ref{alg:gen_fast_diagonalization}.

We assume that Step~1 can be computed with FD method.
Notice that, at Step~2 one computes the generalized eigendecomposition of
the pair $(\mbf A_t,\mbf M_t)$ with cost $O(N_t^3)$ and the
diagonal matrix $\mbf T$ with cost $O(N)$. 
It should be mentioned that in certain cases (see e.g. \cite{Buzbee1970, Haidvogel1979, Shen1997, kressner2023improved, Montardini2023, gander2024new}) 
the matrix with $\mbf{U}_{s,l}$, $l=1,\ldots,d$, and/or with $\mbf{U}_{t}$ appearing at steps 3 and 5 can be performed by exploiting the Fast Fourier Transform, resulting in an almost linear computational cost. However, for simplicity here we assume that matrix products are computed in the standard way.
Table~\ref{tab:fast_diagonalization_costs} describes the overall computational complexity for the DT method. 

\begin{table}
    \centering
    \setlength{\tabcolsep}{16pt}  
    \begin{tabular}{cc|cc}
        \svhline\noalign{\smallskip}
        Step~1 &  Step~2 & Step~3,5 & Step~4 \\
        \noalign{\smallskip}\svhline\noalign{\smallskip}
        $\sum_{l=1}^dN_{s,l}^3$ & $N_t^3 + N $ & $N(N_t + \sum_{l=1}^d N_{s,l})$ & $N$\\
\noalign{\smallskip}\svhline\noalign{\smallskip}  
    \end{tabular}
\caption{Order of complexity of the steps of Algorithm~\ref{alg:gen_fast_diagonalization} for the DT method.}
    \label{tab:fast_diagonalization_costs}    
\end{table}

While this method is very effective in some settings, see \cite{maday2008parallelization}, it may in other cases lead to instabilities and inaccurate solutions. In the remainder of this section, we present two distinct time discretizations of \eqref{eq:model_problem} that highlight these situations.
We first consider a time stepping scheme, and then the Galerkin approximation method. 

\subsection{Finite differences}
In \cite{maday2008parallelization}, the authors considered an implicit Euler scheme in time coupled with a finite difference approximation in space.  
For a fixed time partition $\{t_n\}_{n=0}^{N_t}$ of the interval $\td = [0,T]$, denote by $\mbf{u}_n \in \R^{N_s}$ 
the unknown numerical approximation of the solution of \eqref{eq:model_problem} in the $N_{s}$ internal grid points 
at the time $t_n$. Analogously $\mbf{f}_n \in \R^{N_s}$ is the vector 
representing the source term. Let $\mbf{A}_s \in \R^{N_s\times N_s}$ be 
the symmetric positive definite matrix representing $A_s$, then the set of 
fully discrete equations is 

\begin{equation*}
    \frac{\mbf{u}_{n+1}-\mbf{u}_{n}}{\tau_n} + \mbf{A}_s \mbf{u}_{n+1} = \mbf{f}_{n+1},
    \quad n = 0,\dots,N_t-1,
\end{equation*}
where $\tau_n := t_{n+1}-t_{n}$, for $n=0,\dots, N_t-1$. Therefore, equation~\eqref{eq:kronecker-structure-discrete} holds with 
\begin{equation}\label{eq:implicit_euler_stencil}
    \mbf{A}_t := 
    \begin{pmatrix}
        \frac{1}{\tau_0}  & &\\
        -\frac{1}{\tau_1} & \frac{1}{\tau_1}  &                      &  & \\
   & -\frac{1}{\tau_2} & \frac{1}{\tau_2} &  & \\
   &                       &                      &  \ddots                     & \\
   &                       &                      & -\frac{1}{\tau_{N_t-1}} & \frac{1}{\tau_{N_t-1}} 
    \end{pmatrix},
\end{equation}
and $\mbf{M}_t$, $\mbf{M}_s$ the identity matrices. 

The matrix $\mathbf{A}_t$ in \eqref{eq:implicit_euler_stencil} is lower triangular and its eigenvalues are its diagonal entries. Hence, there exist an invertible matrix $\mbf{U}_t \in \C^{N_t \times N_t}$ that satisfies \eqref{eq:eigen-decomposition-t} \emph{only if} $\tau_i \neq \tau_j$ for every $i\neq j$. 
To ensure this condition, Maday and Rønquist \cite{maday2008parallelization} proposed the use of a geometric partition of the kind 
$$
\tau_n = \beta^{n-1}\tau_1, \quad n \ge 1,
$$
where $\beta>1$ is a free parameter and $\tau_1$ is the first time step. 

By using instead a time step closer to the uniform one, it leads indeed to bad conditioning and instability. The study of these phenomena, in conjunction with the approximation properties of the scheme, has been explored in the work \cite{gander2016direct}.

\subsection{Galerkin discretization}\label{subseq:galerkin-example}
Consider $\sd = [0,1]^{d}$ and let 
$V_s :=  \spann{(B_{i}| i=1,\dots,N_s)}$ be the discrete space of tensor product B-splines of degree $p_s$ and $C^{p_s-1}$ regularity. 
Analogously, $\td = [0,1]$ is the time interval and we denote by $V_t = \spann{(b_{j}| j=1,\dots,N_t)}$ the spline space of degree $p_t$ and regularity $C^{p_t-1}$. For details on spline functions we refer the reader to \cite{deboor}.% We assume maximum smoothness, that is  and  regularity of the basis functions.%, with mesh size $h$. % defined in \cite[Section 2.2]{loli2020efficient},
We assume that  $f \in L^2(\td; H^{-1}(\sd))$ and consider the well posed space-time Galerkin discretization method for \eqref{eq:model_problem} given in \cite[Section 3.2]{loli2020efficient} with the discrete space $V=V_s \otimes V_t$. 
\begin{figure}[H]
    \centering
    \includegraphics[width = 0.85\linewidth]{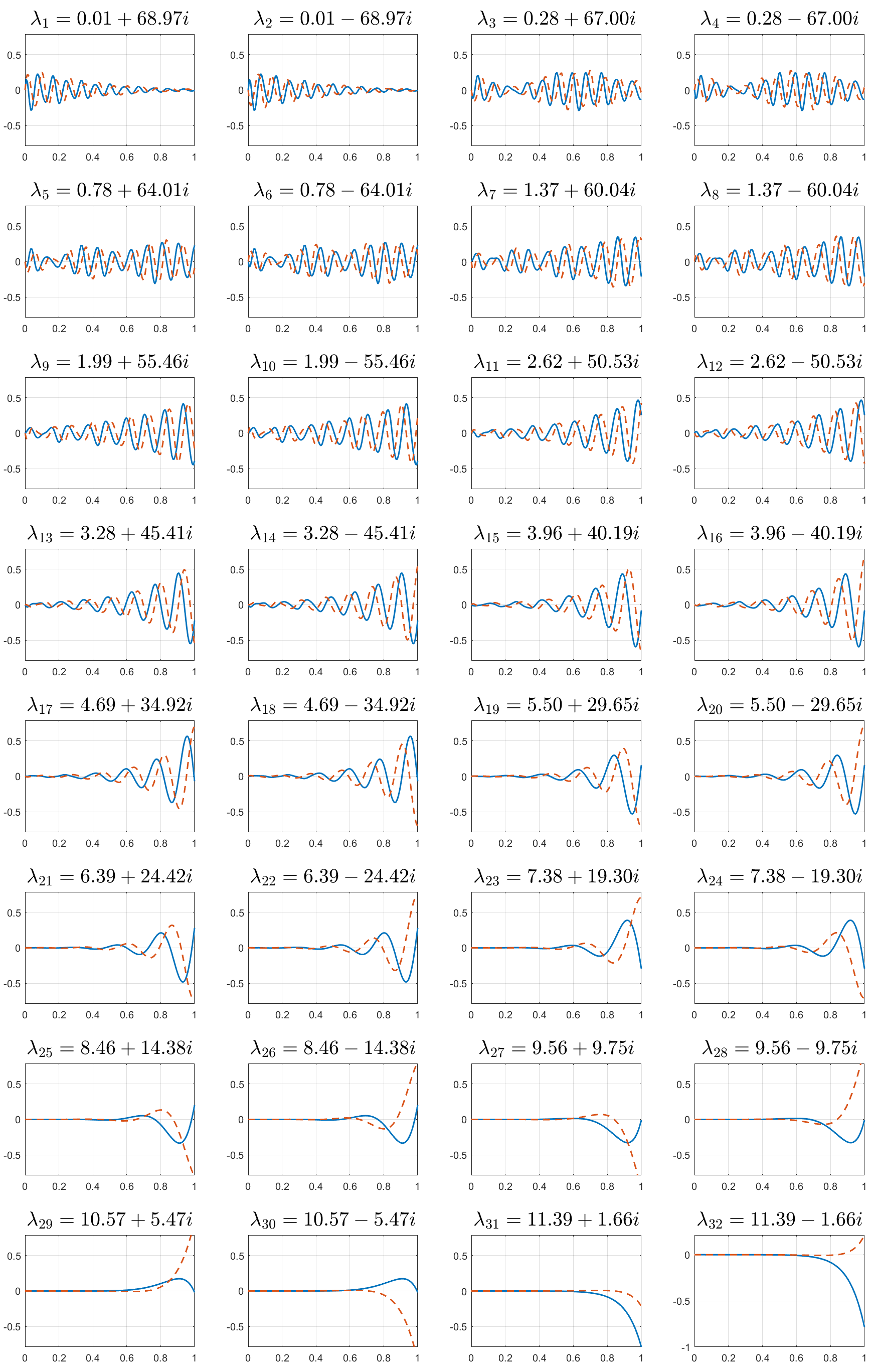}
    \caption{Generalized eigenvectors for the pencil $(\mbf{A}_t, \mbf{M}_t )$, with associated eigenvalues for $p_t =3$ and $N_{t}=32$. The real part is in solid line, while the imaginary part is in dashed line.}
    \label{fig:autovettori_instabili}
\end{figure}
The arising linear system has the structure highlighted in \eqref{eq:kronecker-structure-discrete}, that is 
\[
    \mbf{A} = \mbf{A}_t\otimes \mbf{M}_s + \mbf{M}_t \otimes \mbf{A}_s,
\]
where space matrices are, for $i,j=1,\dots,N_s $,
\begin{equation}
    \label{heat_eq:space_mat}
    [\mbf{A}_s]_{i,j}  =  \int_{\sd} \nabla  B_{j}(\vect{x})\cdot \nabla  B_{i}(\vect{x}) \ \d\sd, 
    \quad  
    [ \mbf{M}_s]_{i,j}  =  \int_{\sd}  B_{j}(\vect{x}) \  B_{i}(\vect{x}) \ \d\sd,
\end{equation}
while the time matrices are for $i,j = 1, \dots, N_t$ 
\begin{equation*}
    [\mbf{A}_t]_{i,j} := \int_{\td} b_{j}'(t)\ b_{i}(t) \dt , \quad [\mbf{M}_t]_{i,j} := \int_{\td} b_{j}(t)\ b_{i}(t) \dt,
\end{equation*}
and are sparse banded with bandwidth $\band =2 p_t+1 \geq1$.% that depends on $p_t$.  

Also in this case, the generalized eigendecomposition of $(\mbf{A}_t,\mbf{M}_t)$ leads to numerical instabilities \cite{loli2020efficient}.  
The eigendecomposition~\eqref{eq:eigen-decomposition-t} is unstable under mesh refinements and degree elevations. 
In particular, a numerical computation reveals that the eigenvectors are far from $\mbf{M}_t$--orthogonality, see Figure~\ref{fig:autovettori_instabili}, for $T=1$, $p_t=3$ and $N_{t}=32$.
In Table~\ref{tab:cond_number}, we report $\kappa_2(\mbf{U}_t)$, varying the polynomial degree $p_t$ and the number of degrees of freedom $N_t$, showing that $\kappa_2(\mbf{U}_t)$ is unstable and
grows exponentially with respect to the spline degree.  
\begin{table}[H]\scriptsize 
	\centering
    \setlength{\tabcolsep}{3pt}
    \begin{tabular}{l!{\vrule width 1pt}*{7}{c|}c}           
    \svhline
    \hbox{\diagbox[linewidth=2\arrayrulewidth]{$N_t$}{$p_t$}
    %\hskip0.8ex\vrule\null\vrule\hskip-0.01ex\vrule
    } & $1$ &  $2$  & $3$ & $4\phantom{^0}$ & $5\phantom{^0}$ & $6\phantom{^0}$ & $7\phantom{^0}$ & $8\phantom{^0}$ \\     
%    $N_{t}$ & $p_t=1$ & $p_t=2$  & $p_t=3$ & $p_t=4$ & $p_t=5$ & $p_t=6$ & $p_t=7$ & $p_t=8$ \\     
    \svhline\noalign{\smallskip}            
    $32$    & $ 1.0 \cdot 10^2 $ & $ 9.8 \cdot 10^2 $ & $ 2.7 \cdot 10^4 $ & $ 3.9 \cdot 10^{4\phantom{0}}  $ & $ 1.5 \cdot 10^{5\phantom{0}}  $ & $ 7.6 \cdot 10^{5\phantom{0}}  $ & $ 4.4 \cdot 10^{6\phantom{0}}  $ & $ 3.7 \cdot 10^{7\phantom{0}}$  \\     
    %\hline   
    $64$    & $ 2.9 \cdot 10^2 $ & $ 4.9 \cdot 10^3 $ & $ 2.8 \cdot 10^5 $ & $ 4.6 \cdot 10^{5\phantom{0}}  $ & $ 3.4 \cdot 10^{6\phantom{0}}  $ & $ 3.6 \cdot 10^{7\phantom{0}}  $ & $ 6.3 \cdot 10^{8\phantom{0}}  $ & $ 2.8 \cdot 10^{9\phantom{0}}$ \\      
    %\hline   
    $128$   & $ 8.7 \cdot 10^2 $ & $ 2.6 \cdot 10^4 $ & $ 1.3 \cdot 10^6 $ & $ 6.1 \cdot 10^{6\phantom{0}}  $ & $ 7.9 \cdot 10^{7\phantom{0}}  $ & $ 1.6 \cdot 10^{9\phantom{0}}  $ & $ 1.6 \cdot 10^{11} $ & $ 2.5 \cdot 10^{11}$ \\ 
    %\hline   
    $256$   & $ 2.6 \cdot 10^3 $ & $ 1.4 \cdot 10^5 $ & $ 1.1 \cdot 10^7 $ & $ 8.5 \cdot 10^{7\phantom{0}}  $ & $ 1.9 \cdot 10^{9\phantom{0}}  $ & $ 7.0 \cdot 10^{10} $ & $ 4.8 \cdot 10^{12} $ & $ 1.6 \cdot 10^{13}$ \\
    %\hline   
    $512$   & $ 8.0 \cdot 10^3 $ & $ 8.0 \cdot 10^5 $ & $ 1.0 \cdot 10^8 $ & $ 1.3 \cdot 10^{9\phantom{0}}  $ & $ 4.9 \cdot 10^{10} $ & $ 3.6 \cdot 10^{12} $ & $ 1.1 \cdot 10^{13} $ & $ 2.6 \cdot 10^{12}$ \\
    %\hline   
    $1024$  & $ 2.6 \cdot 10^4 $ & $ 4.7 \cdot 10^6 $ & $ 9.6 \cdot 10^8 $ & $ 2.2 \cdot 10^{10} $ & $ 1.3 \cdot 10^{12} $ & $ 1.6 \cdot 10^{13} $ & $ 4.2 \cdot 10^{13} $ & $ 7.3 \cdot 10^{12}$ \\
    \noalign{\smallskip}\svhline\noalign{\smallskip}
\end{tabular}
    \caption{$\kappa_2 (\mbf{U}_t)$ in the generalized eigendecomposition of $(\mbf{A}_t,\mbf{M}_t)$ for different degrees $p_t$ and sizes $N_{t}$.}                          
    \label{tab:cond_number}       
\end{table}

\section{Stable time factorizations}\label{sec:stable_time_factorizations}
This section describes the three alternatives to the DT method, denoted LU, AR and LR, considered in Algorithm~\ref{alg:gen_fast_diagonalization}.
Recall that the generalized eigendecomposition \eqref{eq:eigen-decomposition-s}
can be computed and the matrix $\mbf{A}$ factorizes as in \eqref{eq:generic-factorization}, that is 
\[
    \mbf{A} = (\widetilde{\mbf{U}}_t\otimes \widetilde{\mbf{U}}_s)^{-1} \mbf{T}(\mbf{U}_t\otimes\mbf{U}_s)^{-1}.
\]
Let us recall that Kronecker product is not commutative, but it holds for instance
$$
\mbf{U}_s \otimes \mbf{U}_t = \mbf{S}_{N_t,N_s} (\mbf{U}_t \otimes \mbf{U}_s) \mbf{S}_{N_t,N_s}^{\top}, 
$$ 
with $\mbf{S}_{N_t,N_s} \in \R^{N \times N}$ a \textit{perfect shuffle} matrix \cite{van2000ubiquitous}. 
Applying the shuffle matrix corresponds to swapping the positions of time and space in the tensor product structure. 
Clearly $\mbf{S}_{N_t,N_s}^{-1} = \mbf{S}_{N_t,N_s}^{\top}$ and denoting by 
$ \over{T} := \mbf{S}_{N_t,N_s} \mbf{T} \mbf{S}_{N_t,N_s}^{\top} $ the commuted matrix, it holds 
\begin{equation}
    \begin{aligned}
        \over{T} &= \mbf{S}_{N_t,N_s} \mbf{T} \mbf{S}_{N_t,N_s}^{\top} \\
                 &= \mbf{S}_{N_t,N_s} 
                        (\widetilde{\mbf{U}}_t \mbf{A}_t \mbf{U}_t \otimes \mbf{I}_{N_s}
                        + \widetilde{\mbf{U}}_t \mbf{M}_t \mbf{U}_t \otimes \mbf{\Lambda}_s ) \mbf{S}_{N_t,N_s}^{\top}  \\
                 &=  \mbf{I}_{N_s} \otimes \widetilde{\mbf{U}}_t \mbf{A}_t \mbf{U}_t
                 + \mbf{\Lambda}_s \otimes \widetilde{\mbf{U}}_t \mbf{M}_t \mbf{U}_t.      
    \end{aligned}        
\end{equation}
The matrix $\over{T}$ has a block diagonal structure, where each block 
$$\over{T}_t (\lambda_i) := \widetilde{\mbf{U}}_t\mbf{A}_t\mbf{U}_t + \lambda_i \widetilde{\mbf{U}}_t\mbf{M}_t\mbf{U}_t$$ 
is associated to an eigenvalue $\lambda_i := [\mbf{\Lambda}_s]_{i,i}$, for $i=1,\dots, N_s$, that is 
\begin{equation}\label{eq:block_diagonal_structure}
    \over{T} = 
    \begin{pmatrix}
        \over{T}_t(\lambda_1) & & &\\
         & \over{T}_t (\lambda_2) & & \\
         & & \ddots & \\
         & & & \over{T}_t(\lambda_{N_s})
    \end{pmatrix},
\end{equation}
and the expression of $\mbf U_t$, $\widetilde{\mbf{U}}_t$ and $\over{T}_t(\lambda)$ depends on the chosen method. 
Moreover, computing $\mbf{T}^{-1}\mbf{y}$ in Step~4 of Algorithm~\ref{alg:gen_fast_diagonalization} is equivalent to compute 
\[
    \mbf{T}^{-1} = \mbf{S}_{N_t,N_s}^{\top} \over{T}\,^{-1} \mbf{S}_{N_t,N_s},
\]
that is to apply twice a shuffle matrix to a vector, and invert a block diagonal matrix $\over{T}$,
that is, to solve $N_s$ independent subproblems. For this reason all the methods allow for a space parallel implementation. Time parallelization is achieved when each subproblem allows for additional parallelism.

When reporting the cost of the different methods, it is assumed 
the generalized eigendecomposition \eqref{eq:eigen-decomposition-s} can be computed with the same cost as in the FD method, sufficient conditions
are \eqref{eq:Ms-tp} and \eqref{eq:kronecker_structure_of_stiffness_matrix}.

\subsection{Banded matrices and LU~factorization}\label{sub_sec:banded}
This approach is well suited if the time matrices $\mbf{A}_t$ and $\mbf{M}_t$ have small bandwidth $\band$. %  and it extends to a broader class of time operators but is not suited for time parallelization. 
In this method
\[
\mbf{U}_t=\widetilde{\mbf{U}}_t=\mbf{I}_{N_t},
\]
so that 
\begin{equation}\label{eq:Tlambda-LU}
    \over{T}_t(\lambda):= \mbf{A}_t+ \lambda\mbf{M}_t.
\end{equation}
The LU-factorization of $\over{T}_t(\lambda)$ preserves the bandwidth
$$
\over{T}_t(\lambda) = 
\begin{pmatrix}
    l_{1,1}(\lambda) &       & & \\
    \vdots  & \ddots& & \\
    l_{\band,1}(\lambda) &       & & \\
            & \ddots& \quad \quad & l_{N_t,N_t}(\lambda)
\end{pmatrix}
\begin{pmatrix}
    u_{1,1}(\lambda)  &   \dots & u_{1,\band}(\lambda) &   \\
             &  \ddots &         & \ddots\\
             &         &         &    \\
             &         &         &   u_{N_t,N_t}(\lambda)
\end{pmatrix},
$$
and it is computed in Step~2 of Algorithm \ref{alg:gen_fast_diagonalization} for each $\lambda_i$, $i=1,\dots,N_s$.
Then the inversion of $\mbf{T}$ in Step~4 of Algorithm~\ref{alg:gen_fast_diagonalization} corresponds to the solution of $N_s$ independent LU-factorized linear systems with bandwidth $\band$.
The overall computational cost of this approach is described in Table~\ref{tab:lu_factorization_costs}.
\begin{table}
    \centering
    \setlength{\tabcolsep}{16pt}  
    \begin{tabular}{cc|cc}
        \svhline\noalign{\smallskip}
        Step~1 & Step~2& Step~3,5 & Step~4 \\
        \noalign{\smallskip}\svhline\noalign{\smallskip}
        $\sum_{l=1}^dN_{s,l}^3$ & $\band^2 N $ & $N \sum_{l=1}^d N_{s,l}$ & $\band   N$ \\
\noalign{\smallskip}\svhline\noalign{\smallskip}  
    \end{tabular}
    \caption{Order of complexity of the steps of Algorithm~\ref{alg:gen_fast_diagonalization} for  the LU method.}
    \label{tab:lu_factorization_costs}    
\end{table}

\subsection{Arrowhead factorization}\label{sub_sub_sec:arrow}

The arrowhead factorization was introduced in \cite{loli2020efficient}. %to solve the problem described in Subsection~\ref{subseq:galerkin-example} where the generalized diagonalization of $\mbf{A}_t$, $\mbf{M}_t$ is numerically unstable. In that case, or 
For any basis such that $b_{j}(T)=0$, $j=1,\dots,N_t-1$, $b_{N_t}(T)\ne 0$, the Galerkin time derivative $\mbf{A}_t$ and the Gram matrix $\mbf{M}_t$ have the block structure:
\begin{equation*}\label{eq:split}
	\mbf{A}_t = 
	\begin{bmatrix}
		\mathring{\mbf{A}}_t & \mbf{a}\\[2pt]
		-\mbf{a}^{\top} & \alpha
	\end{bmatrix}
	\quad \text{and} \quad
	\mbf{M}_t=
	\begin{bmatrix}
		\mathring{\mbf{M}}_t & \mbf{m}\\[2pt]
		\mbf{m}^{\top} & \mu
	\end{bmatrix},
\end{equation*}
where $\mathring{\mbf{M}}_t$ is symmetric, $\mathring{\mbf{A}}_t$ is skew-symmetric and 
\[
\alpha = \int_{0}^{T} b_{N_t}'\ b_{N_t} \dt  =  \frac{b_{N_t}(T)^2}2>0.
\]
Therefore, there exists $\mathring{\mbf{U}}_t$ satisfying
\begin{equation*}
	\mathring{\mbf{U}}_t^*  \mathring{\mbf{A}}_t \mathring{\mbf{U}}_t = \mathring{\mbf{D}}_t  
	\quad \text{ and } \quad 
	\mbf{U}_t^*  \mbf{M}_t \mbf{U}_t = \mbf{I}_{N_t}.
\end{equation*}
Let us define $\mbf{U}_t$ as
\begin{equation}\label{eq:U_t}
	\mbf{U}_t := 
	\begin{bmatrix}
		\mathring{\mbf{U}}_t & -\rho \mathring{\mbf{M}}^{-1}_t \mbf{m} \\
		\mbf{0}^\top & \rho
	\end{bmatrix},
\end{equation}
with $\mbf 0$ being the null vector in $\R^{N_t-1}$ and
\begin{equation*}
	\rho := \left( \mu - \mbf{m}^{\top} \mathring{\mbf{M}}^{-1}_t \mbf{m} \right)^{-1/2} = \left( \| b_{N_t} \|^2_{L^2(0,T)} - \| \Pi (b_{N_t}) \|^2_{L^2(0,T)} \right)^{-1/2},
\end{equation*}
where 
\begin{equation*}
	\Pi: L^2(0,T) \rightarrow \spann\{ b_i : i=1,\ldots,N_t-1 \}
\end{equation*}
denotes the $L^2$-projection. Since $b_{N_t} \notin \spann\{ b_i : i=1,\ldots,N_t-1 \}$, we highlight that
\begin{equation*}
	\mu - \mbf{m}^{\top} \mathring{\mbf{M}}^{-1}_t \mbf{m} > 0.
\end{equation*}
Choosing $\mbf{U}_t$ as in \eqref{eq:U_t}, we obtain:
\begin{equation*}
	\widetilde{\mbf{U}}_t = \mbf{U}_t^*, \quad \mbf{U}_t^* \mbf{A}_t \mbf{U}_t = \mbf{\Delta}_t \quad \text{and} \quad \mbf{U}_t^* \mbf{M}_t \mbf{U}_t = \mbf{I}_{N_t},
\end{equation*}
where $\mbf{\Delta}_t$ is a complex arrowhead matrix, characterized by non-zero entries only on the diagonal, in the last row, and in the last column, resembling an arrow pointing down-right. Specifically, the expression for $\mbf{\Delta}_t$ is
\begin{equation*}\label{eq:time_eig_2}
	\mbf{\Delta}_t:=\mbf{U}_t^*\mbf{A}_t \mbf{U}_t   =  \begin{bmatrix}
		\mathring{\mbf{D}}_t& \mbf{g}\\[2pt]
		-\mbf{g}^* 
		& \sigma
	\end{bmatrix},
\end{equation*}
where 
\[\mbf{g}:=\mathring{{\mbf{U}}}_t^*\begin{bmatrix} \mathring{\mbf{A}}_t & \mbf{a} \end{bmatrix} \begin{bmatrix}
	-\rho \mathring{\mbf{M}}^{-1}_t \mbf{m} \\ \rho
\end{bmatrix},\qquad\text{and}\qquad\sigma:=\begin{bmatrix} -\rho \mbf{m}^{\top} \mathring{\mbf{M}}^{-1}_t & {\rho} \end{bmatrix} \mbf{A}_t \begin{bmatrix}
	-\rho \mathring{\mbf{M}}^{-1}_t \mbf{m} \\ \rho
\end{bmatrix}.\]

With this method the blocks of $\over{T}$ in \eqref{eq:block_diagonal_structure} 
have the same arrowhead structure of $\mathbf{\Delta}_t$, and the LU~factorization 
of $\over{T}_t(\lambda)$ is as follows
\begin{equation}\label{eq:stable_arrow_factorization}
    \over{T}_t (\lambda) = \mbf{\Delta}_t + \lambda \mbf{I}_{N_t}
    =\begin{pmatrix}
        1 & & & \\
        & \ddots & & \\
        & & 1 & \\
        -\frac{g_1^*}{\delta_1+\lambda} & \dots & -\frac{g_{N_t-1}^*}{\delta_{N_t-1}+\lambda} & 1
    \end{pmatrix}
    \begin{pmatrix}
        (\delta_1+\lambda) & & & g_1 \\
        & \ddots & & \vdots \\
        & & (\delta_{N_t-1}+\lambda) & g_{N_t-1}\\
        & & & s(\lambda)
    \end{pmatrix}
\end{equation}
where $g_j = [\mbf{g}]_j$, $\delta_j = [\mbf{\Delta}_t]_{j,j}$ for 
$j = 1,\dots, N_t-1$, $\sigma = [\mbf{\Delta}_t]_{N_t,N_t} $ and $$s(\lambda) = (\sigma + \lambda) + \sum_{j=1}^{N_t-1} g_j^* (\delta_j + \lambda)^{-1}g_j.$$ 
The application of the inverse of $\over{T}$ is now parallelizable in space and time. Parallelization in space is of course achievable by inverting each block independently. Parallelization in time is achieved by backward substitution, inverting first the last row of \eqref{eq:stable_arrow_factorization}, and then the previous $N_t-1$ rows concurrently. 

The setup of the AR method incudes the computation of the LU factors for each $\lambda_i \in \mbf{\Lambda}_s$ whose computational cost is $O(N)$.
The inversion of $\mbf{T}$ in Step~4 of Algorithm~\ref{alg:gen_fast_diagonalization}
is again of order $O(N)$ due to the form of the LU factorization, so that the dominant costs are given by the matrix vector multiplications in Step~3 and Step~5.
The overall computational costs are reported in Table~\ref{tab:arrow_factorization_costs}.
\begin{table}
    \centering
    \setlength{\tabcolsep}{16pt}  
    \begin{tabular}{cc|cc}
        \svhline\noalign{\smallskip}
            Step~1 & Step~2 & Step~3,5 & Step~4\\
            \noalign{\smallskip}\svhline\noalign{\smallskip}
            $\sum_{l=1}^dN_{s,l}^3$ & $N_t^3+N $ & $N(N_t + \sum_{l=1}^d N_{s,l})$& $N$   \\
\noalign{\smallskip}\svhline\noalign{\smallskip} 
    \end{tabular}
    \caption{Order of complexity of the steps of Algorithm~\ref{alg:gen_fast_diagonalization} for  the AR method.}
    \label{tab:arrow_factorization_costs} 
\end{table}    

\begin{figure}[H]
    \centering
    \includegraphics[width = 0.85\textwidth]{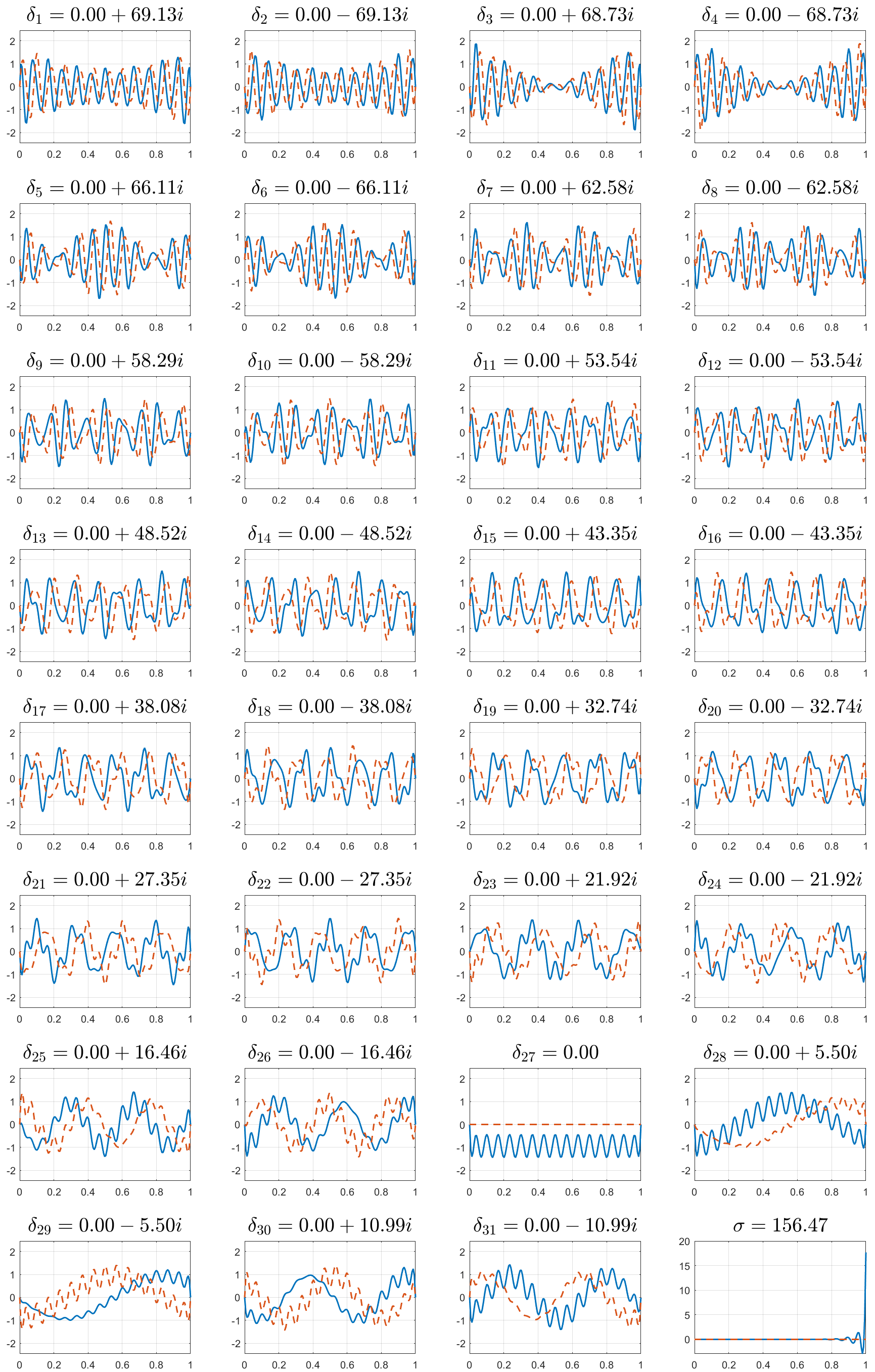}
    \caption{Real part (solid line) and imaginary part (dashed line) of splines corresponding to the columns of $\mbf{U}_t$, as function of time $t \in [0,1]$, with associated diagonal entry in $\mbf{\Delta}_t$. Discretization with $p_t = 3$ and $N_{t} = 32$.}
    \label{fig:autovettori_stabili}
\end{figure}  

  \begin{table}[H]\scriptsize   
  	\centering                    
    \setlength{\tabcolsep}{3pt}
    \begin{tabular}{l!{\vrule width 1pt}*{7}{c|}c}           
        \svhline
        \hbox{\diagbox[linewidth=2\arrayrulewidth]{$N_t$}{$p_t$}
        %\hskip0.8ex\vrule\null\vrule\hskip-0.01ex\vrule
        } & $1$ &  $2$  & $3$ & $4$ & $5$ & $6$ & $7$ & $8$ \\ 
        \svhline\noalign{\smallskip}            
        32    & $2.0 \cdot 10^0 $ & $3.3 \cdot 10^0 $ & $ 5.2 \cdot 10^0$ & $ 8.3 \cdot 10^0$ & $ 1.3 \cdot 10^1$ & $ 2.2  \cdot 10^1$ & $ 3.6  \cdot 10^1$ & $ 5.9  \cdot 10^1$ \\   
 		%\hline    
 		64    & $2.0 \cdot 10^0 $ & $3.3 \cdot 10^0 $ & $ 5.2 \cdot 10^0$ & $ 8.3 \cdot 10^0$ & $ 1.3 \cdot 10^1$ & $ 2.2  \cdot 10^1$ & $ 3.6  \cdot 10^1$ & $ 5.9  \cdot 10^1$ \\   
 		%\hline    
 		128   & $2.0 \cdot 10^0 $ & $3.3 \cdot 10^0 $ & $ 5.2 \cdot 10^0$ & $ 8.3 \cdot 10^0$ & $ 1.3 \cdot 10^1$ & $ 2.2  \cdot 10^1$ & $ 3.6  \cdot 10^1$ & $ 5.9  \cdot 10^1$ \\   
 		%\hline    
 		256   & $2.0 \cdot 10^0 $ & $3.3 \cdot 10^0 $ & $ 5.2 \cdot 10^0$ & $ 8.3 \cdot 10^0$ & $ 1.3 \cdot 10^1$ & $ 2.2  \cdot 10^1$ & $ 3.6  \cdot 10^1$ & $ 5.9  \cdot 10^1$ \\  
 		%\hline    
 		512   & $2.0 \cdot 10^0 $ & $3.3 \cdot 10^0 $ & $ 5.2 \cdot 10^0$ & $ 8.3 \cdot 10^0$ & $ 1.3 \cdot 10^1$ & $ 2.2  \cdot 10^1$ & $ 3.6  \cdot 10^1$ & $ 5.9  \cdot 10^1$ \\   
 		%\hline    
 		1024  & $2.0 \cdot 10^0 $ & $3.3 \cdot 10^0 $ & $ 5.2 \cdot 10^0$ & $ 8.3 \cdot 10^0$ & $ 1.3 \cdot 10^1$ & $ 2.2  \cdot 10^1$ & $ 3.6  \cdot 10^1$ & $ 5.9  \cdot 10^1$ \\ 
 		\noalign{\smallskip}\svhline\noalign{\smallskip}
 	\end{tabular}
    \caption{$\kappa_2 (\mbf{U}_t)$ in arrowhead factorization for different degrees $p_t$ and sizes $N_{t}$.}
    \label{tab:cond_number_time_new}
 \end{table}

Finally, going back to the example of Subsection~\ref{subseq:galerkin-example} where $\mbf{M}_t$ is the Gram matrix, it is known that 
$\kappa_2(\mbf{M}_t)$ is bounded uniformly with respect to the spline knot vector, see for instance~\cite[Lemma 1]{loli2022easy}. %\cite{LYCHE1978202,SCHERER1999217}.
By construction, the columns of $\mbf{U}_t$ are $\mbf{M}_t$--orthonormal, 
which implies that $\kappa_2\left(\mbf{U}_t\right)=\sqrt{\kappa_2(\mbf{M}_t)}$. 
This was numerically confirmed in \cite{loli2020efficient} and the computed $\kappa_2(\mbf{U}_t)$ %for $T=1$, different spline degrees $p_t$ and different numbers of degrees of freedom $N_{t}$ 
are reported in Table~\ref{tab:cond_number_time_new}. Note that, while the conditioning grows with respect to the degree,
it remains reasonably small for all the degrees of interest. 
Finally, Figure~\ref{fig:autovettori_stabili} shows the plot of spline functions corresponding to the columns of 
$\mbf{U}_t$ with associated diagonal entries $\delta_i = [\mbf{\Delta}_t]_{i,i}$, $i = 1,\dots, N_t-1$ 
and  $\sigma = [\mbf{\Delta}_t]_{N_t,N_t} $, for $p_t = 3$ and uniform partition with $N_{t}=32$.

\subsection{Low-rank modification}\label{sub_sub_sec:low_rank}
The arrowhead factorization is a way to cope with the
lack of skew-symmetry in $\mbf{A}_t$ due to a single positive coefficient in 
the diagonal entries.
In that case a skew-symmetric matrix is obtained from $\mbf{A}_t$ by
removing the positive entry, i.e., by a rank-1 modification.

A generalization of the arrowhead factorization that allows for broader set of pencils $(\mbf{A}_t,\mbf{M}_t)$ is to split $\mbf{A}_t$ as
\[\mbf{A}_t  = \widetilde{\mbf{A}}_t + \mbf{R}_t,\]
where the matrix $\mbf{R}_t$ has rank-$\rank$, $1 \le \rank \le N_t$, and the pair $(\widetilde{\mbf{A}}_t, \mbf{M}_t)$ admits the generalized eigendecomposition 
\[ 
        \widetilde{\mbf{A}}_t \mbf{U}_t = \mbf{M}_t \mbf{U}_t \mbf{\Lambda}_t,    
\]
with  $\mbf{U}_t$ an $\mbf{M}_t$-orthonormal matrix \eqref{eq:orthogonality}, or equivalently
\begin{equation} \label{eq:time_eig_lr}
    \widetilde{\mbf{U}}_t = \mbf{U}_t^*, \quad 
    {\mbf{U}}_t^*  \widetilde{\mbf{A}}_t \mbf{U}_t = \mbf{\Lambda}_t  
    \quad \text{ and } \quad 
    {\mbf{U}}_t^*  \mbf{M}_t \mbf{U}_t = \mbf{I}_{N_t}.
\end{equation}

Let $\mbf{R}_t = \mbf{F}^{\top}\mbf{G}$, for some $\mbf{F},\mbf{G} \in \mathbb{R}^{\rank \times N_t}$, then 
the blocks of $\over{T}$ for this method take the form 
\begin{equation}\label{eq:low_rank_perturbed_matrix}
    \over{T}_t(\lambda) = \mbf{\Lambda}_t + \lambda \mbf{I}_{N_t} + \mbf{U}_t^*{\mbf{F}^{\top}\mbf{G}}\mbf{U}_t, %= \mbf{\Lambda}_t + \lambda \mbf{I}_{N_t} + \mbf{U}_t^*{\mbf{F}\mbf{G}}\mbf{U}_t 
\end{equation}
that is, $\over{T}_t(\lambda)$ is a low-rank modification of
a diagonal matrix $\mbf{D}_\lambda:=\mbf{\Lambda}_t + \lambda \mbf{I}_{N_t}$ and 
$\over{T}_t(\lambda)^{-1}$ can be applied
using the Sherman--Morrison--Woodbury formula 
\begin{equation}\label{eq:Sherman-Morrison-Woodbury}
    \left(  \mbf{D}_\lambda +   \mbf{U}_t^*\mbf{F}^{\top}\mbf{G}\mbf{U}_t \right)^{-1} =  
            \mbf{D}_\lambda^{-1} - 
            \mbf{D}_\lambda^{-1}  \mbf{U}_t^*\mbf{F}^{\top}
 \left(\mbf{I}_{\rank} +   \mbf{G}\mbf{U}_t  \mbf{D}_\lambda^{-1}  \mbf{U}_t^* \mbf{F}^{\top} \right)^{-1}
       \mbf{G} \mbf{U}_t   \mbf{D}_\lambda^{-1}.
\end{equation}

Thus, in the setup phase of Algorithm~\ref{alg:gen_fast_diagonalization},
the following small matrices in $\mathbb{R}^{\rank \times \rank}$ are precomputed
for each $\lambda_i \in \mbf{\Lambda_s}$
\[ {\mbf{C}}_\lambda:=\left(\mbf{I}_{\rank} + \mbf{G} \mbf{U}_t \mbf{D}_\lambda^{-1} \mbf{U}_t^* \mbf{F}^{\top}\right)^{-1}.\]
This costs $O(\rank N_t^2)$  to compute $\mbf{U}_t^*\mbf{F}^{\top}$ and $\mbf{G}\mbf{U}_t$
plus $O(\rank^2 N)$ to compute the matrix in parenthesis and $O(N_s\rank^3)$ to invert it.

Applying the Sherman--Morrison--Woodbury formula is done in few steps as described in Algorithm~\ref{alg:low-rank-solve},
while the overall computational cost of the LR method is reported in Table~\ref{tab:Sherman_Morrison_Woodbury_factorization_costs}. 
The steps of Algorithm \ref{alg:low-rank-solve} are sequential but each step is a matrix vector multiplication implying that this method is parallelizable in both time and space.

\begin{algorithm}[H]
    \caption{Sherman--Morrisson--Woodbury.}\label{alg:low-rank-solve}
    \begin{algorithmic}[1]
    \State Compute $\mbf{v}_1  := \mbf{D}_\lambda^{-1} \mbf{y}$ \hfill with cost $O(N_t)$
    \State Compute $\mbf{v}_2  := \mbf{G}\mbf{U}_t\mbf{v}_1$\hfill with cost $O(\rank N_t)$
    \State Compute $\mbf{v}_3  := {\mbf{C}}_\lambda\mbf{v}_2$\hfill with cost $O(\rank^2)$
    \State Compute $\mbf{v}_4  := \mbf{U}_t^*\mbf{F}^{\top}\mbf{v}_3$\hfill with cost $O(\rank N_t)$
    \State Compute $\mbf{v}_5  := \mbf{D}_\lambda^{-1}\mbf{v}_4$\hfill with cost $O( N_t)$
    \State Return  $\mbf{v}_1 + \mbf{v}_5$\hfill with cost $O( N_t)$
    \end{algorithmic}
\end{algorithm}

\begin{table}
	\centering
    \setlength{\tabcolsep}{16pt}
    \begin{tabular}{cc|cc}
        \svhline\noalign{\smallskip}
        Step~1 & Step~2  & Step~3,5 & Step~4\\
        \noalign{\smallskip}\svhline\noalign{\smallskip}
        $\sum_{l=1}^dN_{s,l}^3$ & $  N_t^3+\rank^2N$  & $ N(N_t + \sum_{l=1}^d N_{s,l})$ & $  \rank N $ \\
        \noalign{\smallskip}\svhline\noalign{\smallskip} 
    \end{tabular}
    \caption{Order of complexity of the steps of Algorithm~\ref{alg:gen_fast_diagonalization} for  the LR method.}
    \label{tab:Sherman_Morrison_Woodbury_factorization_costs} 
\end{table}    

Note that, for the LR method there are multiple choices of low-rank modifications $\mbf{R}_t$, and its factorizations $\mbf{R}_t = \mbf{F}^{\top}\mbf{G}$. 
For example, in the discretization of Section~\ref{subseq:galerkin-example}, 
a possible first choice is the rank-1 modification obtained by considering: % only the last diagonal entry of $\mbf{A}_t$, that is 
$$ 
\mbf{R}_t:= 
    \begin{pmatrix}
     \mbf{O} &   \mbf{0}    \\
     \,\mbf{0^{\top}}  &  \alpha \,
    \end{pmatrix},     
    \qquad 
    \mbf{F}^{\top} := \alpha\mbf{e}_{N_t},
    \qquad
    \mbf{G} := \mbf{e}_{N_t}^{\top},
$$
with $\alpha = [\mbf{A}_t]_{N_t,N_t}$, $\mbf{O}\in \R^{(N_t-1)\times( N_t-1)}$ the null matrix,  % Its factorization is explicitly given by the choice $\mbf{F}^{\top} := \alpha\mbf{e}_{N_t} $, and $\mbf{G} := \mbf{e}_{N_t}^{\top}$, with 
and $\mbf{e}_{N_t}$ denoting the last element of the canonical basis of $\R^{N_t}$.
An alternative decomposition could be: 
$$ 
\mbf{R}_t:= 
    \begin{pmatrix}
     \mbf{O}    & \mbf{a} \\
     -\mbf{a}^{\top} & \alpha \,
    \end{pmatrix}, 
    \qquad 
    \mbf{F}^{\top}=\begin{pmatrix}\mbf{a} & \mbf{0}\\\alpha & 1\end{pmatrix}
    \qquad
    \mbf{G}=\begin{pmatrix}\mbf{0}^{\top}&1\\-\mbf{a}^{\top}&0\end{pmatrix}
$$
which is a rank-2 modification.

\section{Numerical tests}\label{sec:heat_equation}
In our numerical tests we compare the computational costs of solving a spline based Galerkin discretization of 
\eqref{eq:model_problem}
by Algorithm~\ref{alg:gen_fast_diagonalization}, 
employing the LU, AR or LR methods.

We consider only sequential executions and we force the use of a single computational thread on a workstation with Intel Core i7-5820K processor, 
running at 3.30 GHz and with 64 GB of RAM. The tests are performed with Matlab R2023a and GeoPDEs toolbox \cite{Vazquez2016}. 

We use the \texttt{decomposition} Matlab function, with \texttt{`banded'} option, to setup LU~factorization 
of the blocks $\over{T}_t(\lambda_i)$. The \texttt{eig} Matlab function is employed to compute the 
generalized eigendecompositions, while the Tensorlab toolbox \cite{Sorber2014} is employed to perform the  multiplications with Kronecker matrices. 
Finally, all the computational times reported on plots and tables, except for Table \ref{tab:preconditioner-performance-times}, have been obtained with \texttt{timeit} built in Matlab function, that performs several runs and saves the average computational time. 

\subsection{Application as a direct solver}\label{subseq:direct_solver_tests}
We considered the Galerkin discretization of Section~\ref{subseq:galerkin-example} with 
the constant source term $f = 1$. For the sake of simplicity, we use splines of the same degree in space and in time  $p_t=p_s=:p$ and uniform knot vectors.

\begin{figure}
    \centering
    \subfloat[][$p = 1$.\label{fig:1}]
      {\includegraphics[ clip=true, scale=0.41]{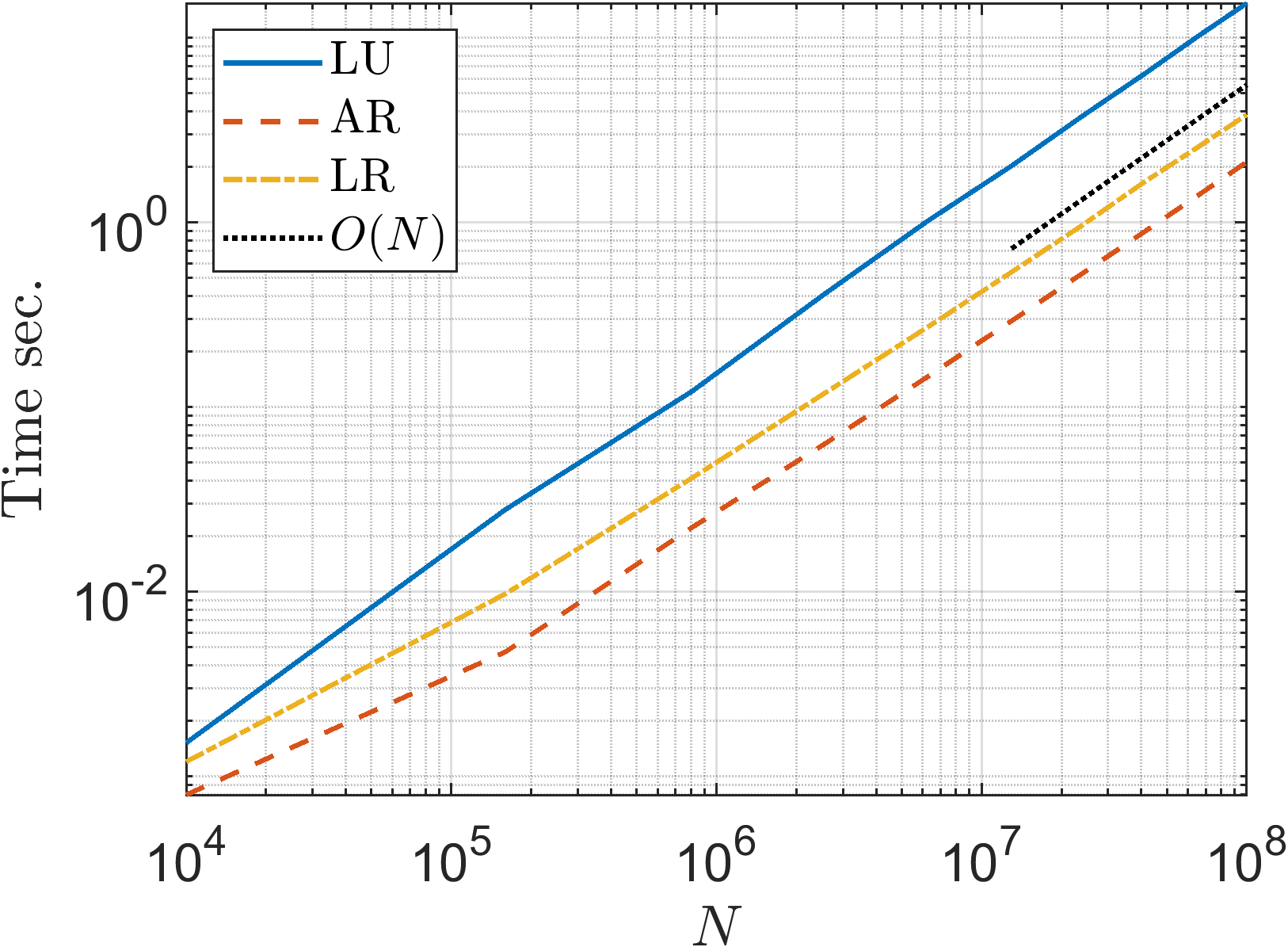}}  \quad 
    \subfloat[][$p = 2$.\label{fig:2}]
      {\includegraphics[ clip=true, scale=0.41]{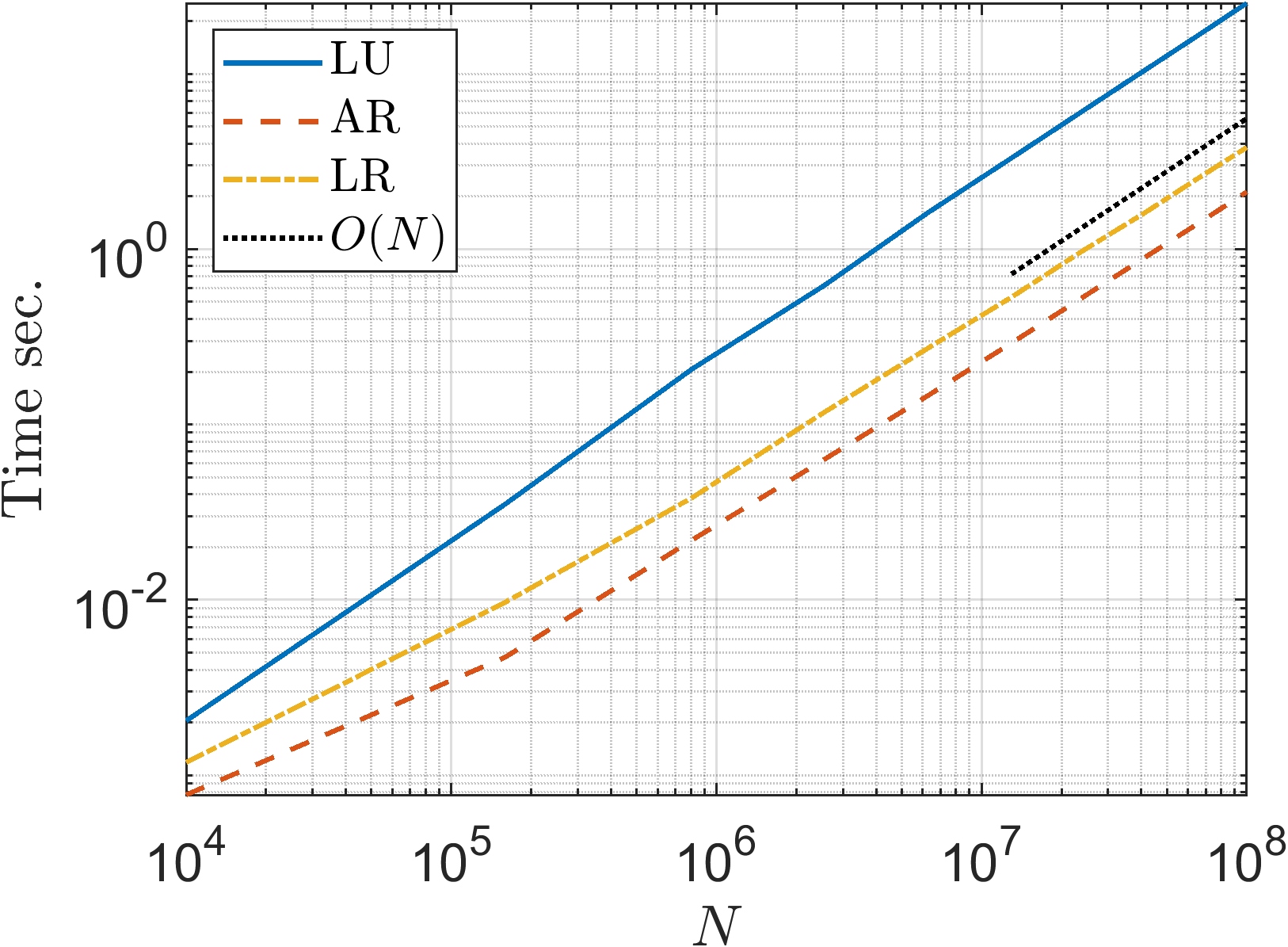}}  \\     
    \subfloat[][$p = 3$.\label{fig:3}]
      {\includegraphics[ clip=true, scale=0.41]{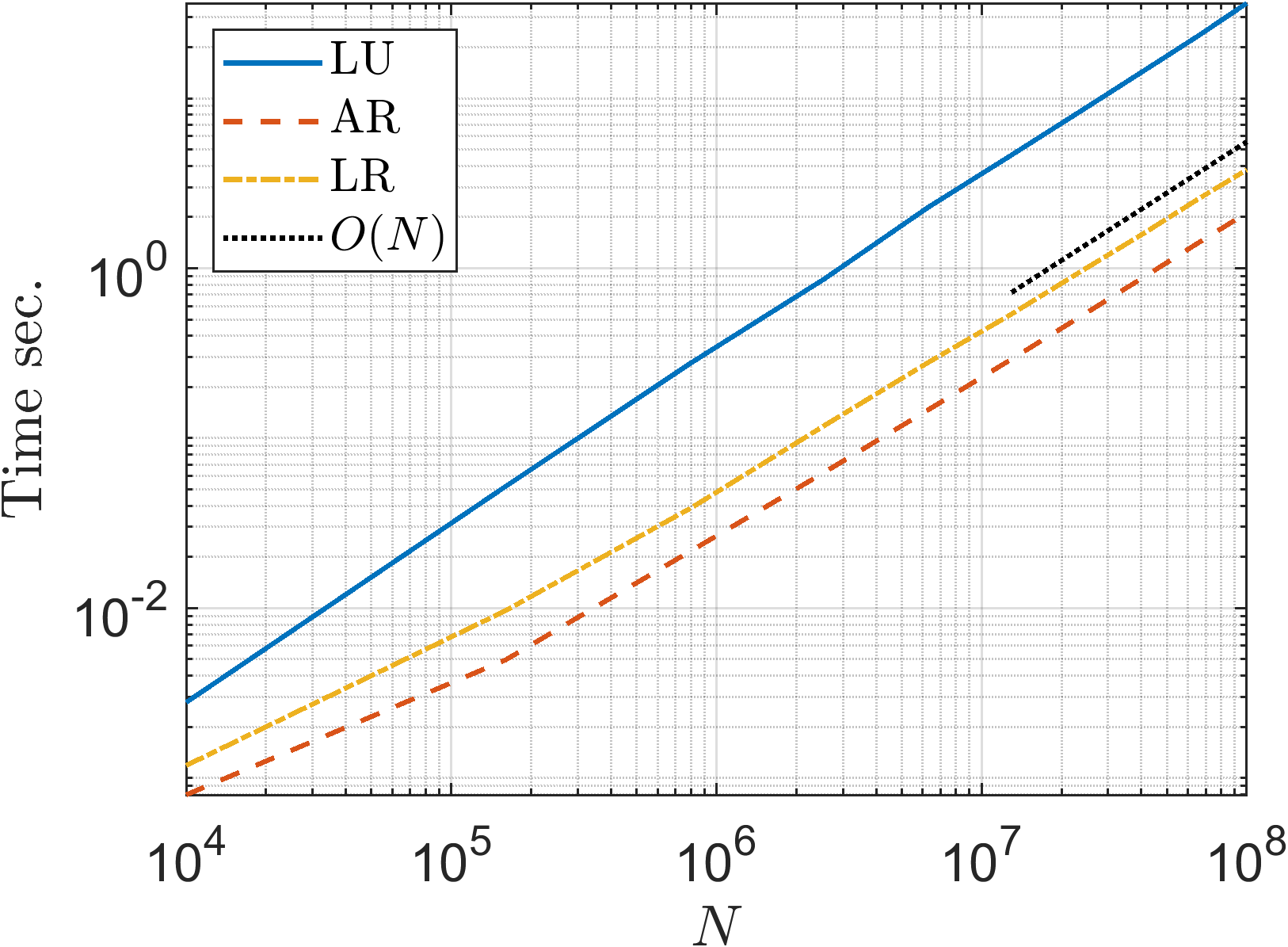}}  \quad 
    \subfloat[][$p = 4$.\label{fig:4}]
      {\includegraphics[ clip=true, scale=0.41]{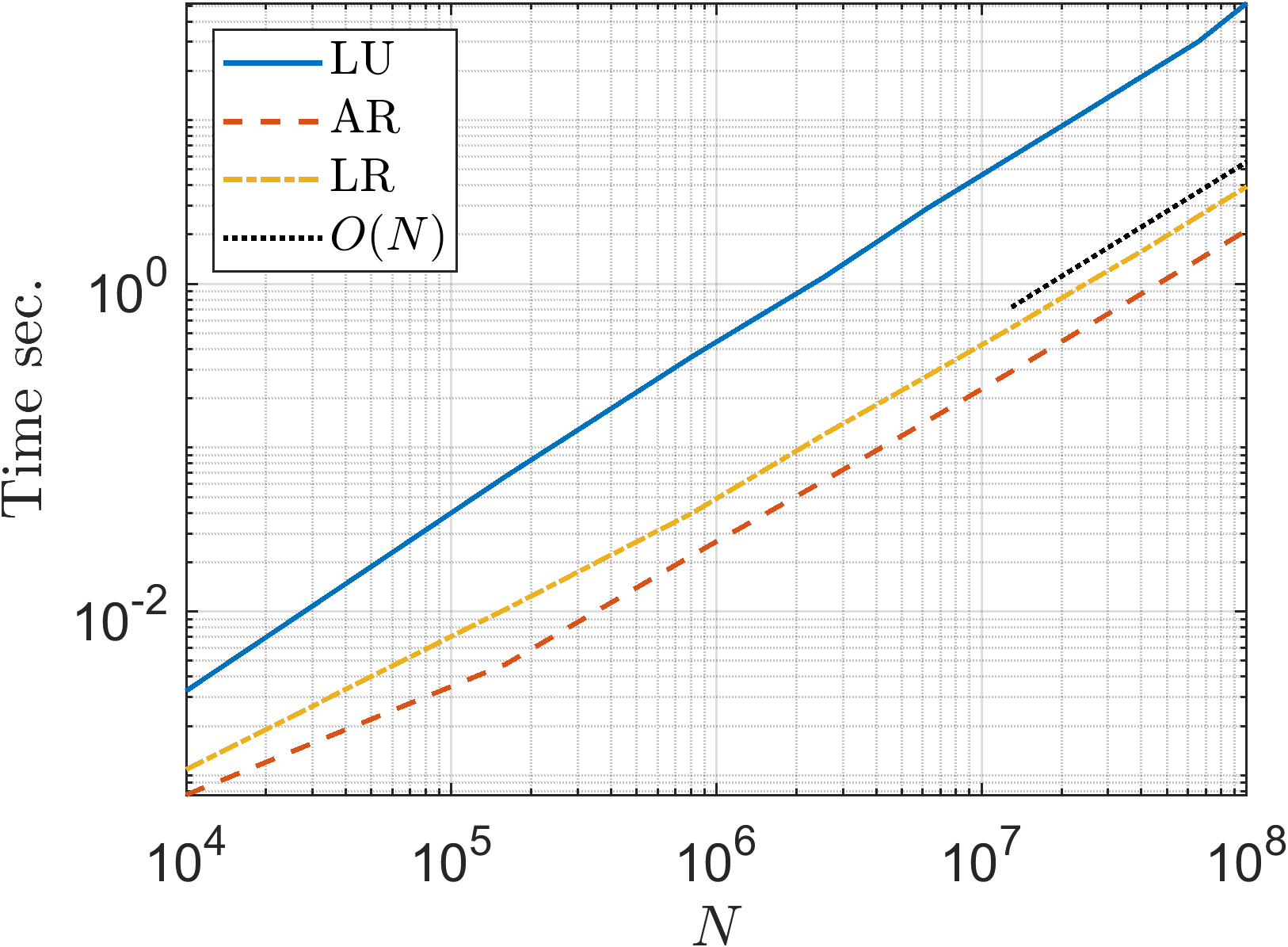}}  %\\
    \caption{Setup cost for $N_s = N_t^3$.   }
    \label{fig:setup_plots}
\end{figure}

\begin{figure}
    \centering
    \subfloat[][$p = 1$.\label{fig:10}]
      {\includegraphics[ clip=true, scale=0.41]{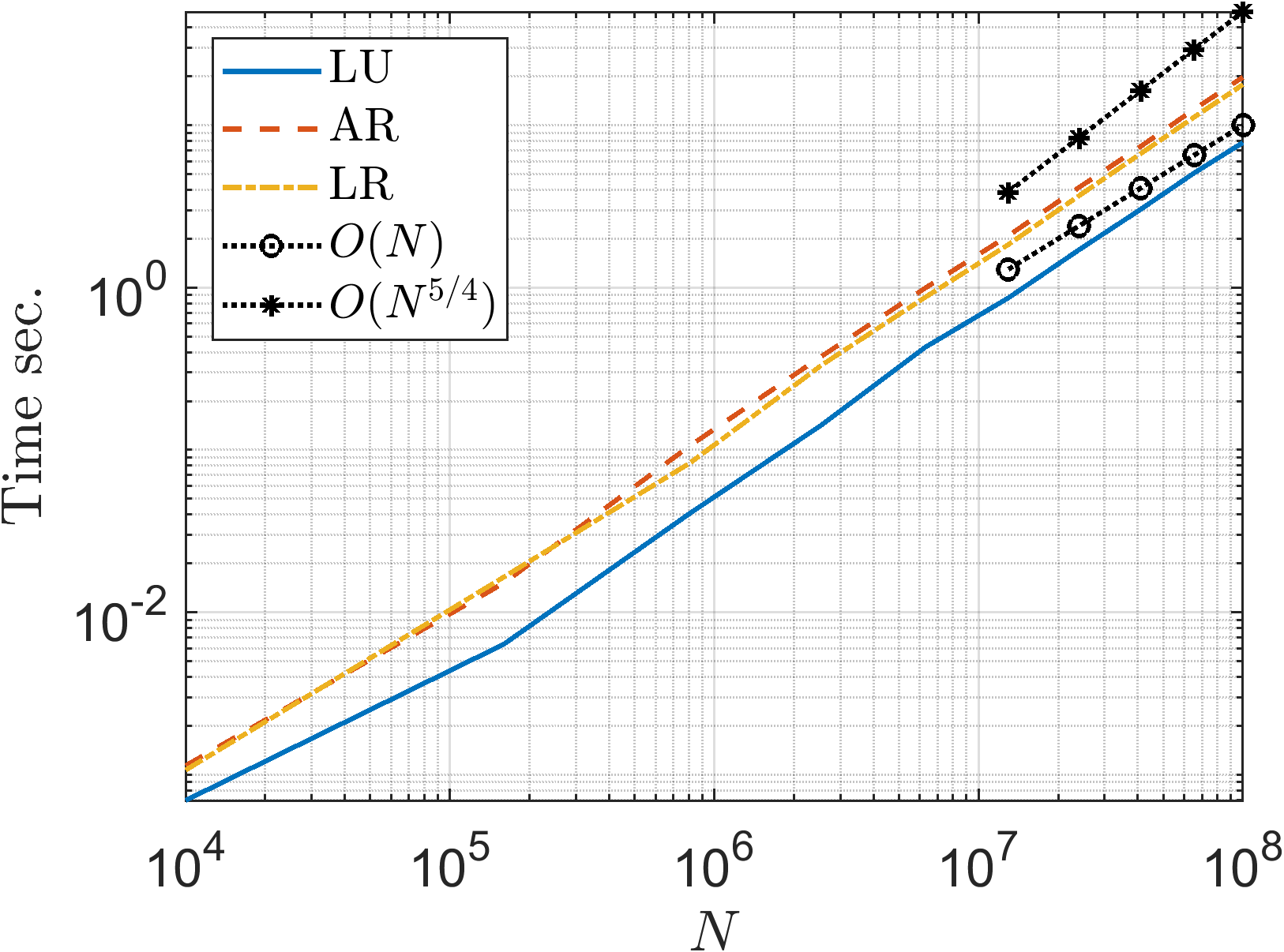}}  \quad 
    \subfloat[][$p = 2$.\label{fig:20}]
      {\includegraphics[ clip=true, scale=0.41]{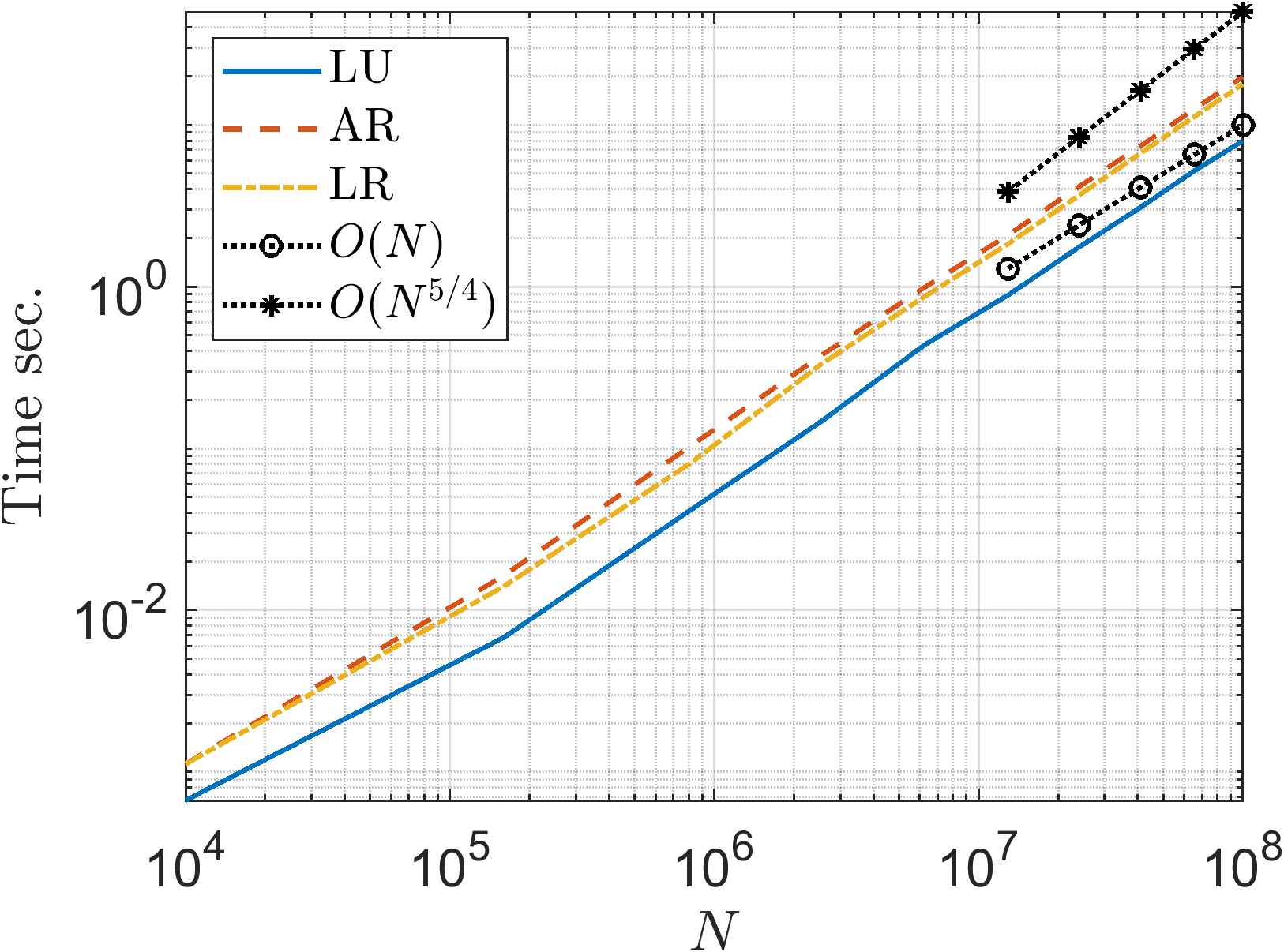}}  \\     
    \subfloat[][$p = 3$.\label{fig:30}]
      {\includegraphics[ clip=true, scale=0.41]{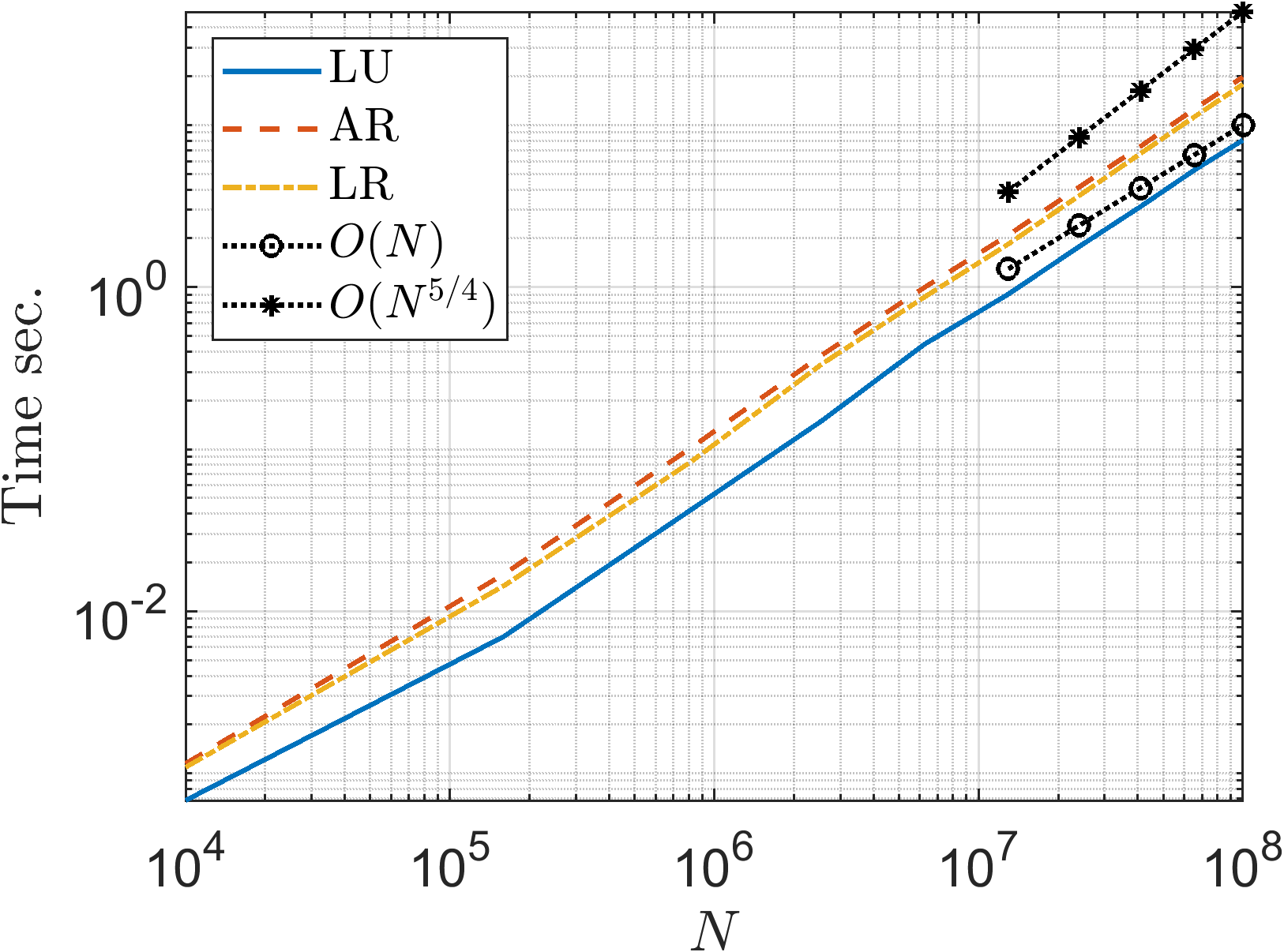}}  \quad 
    \subfloat[][$p = 4$.\label{fig:40}]
      {\includegraphics[ clip=true, scale=0.41]{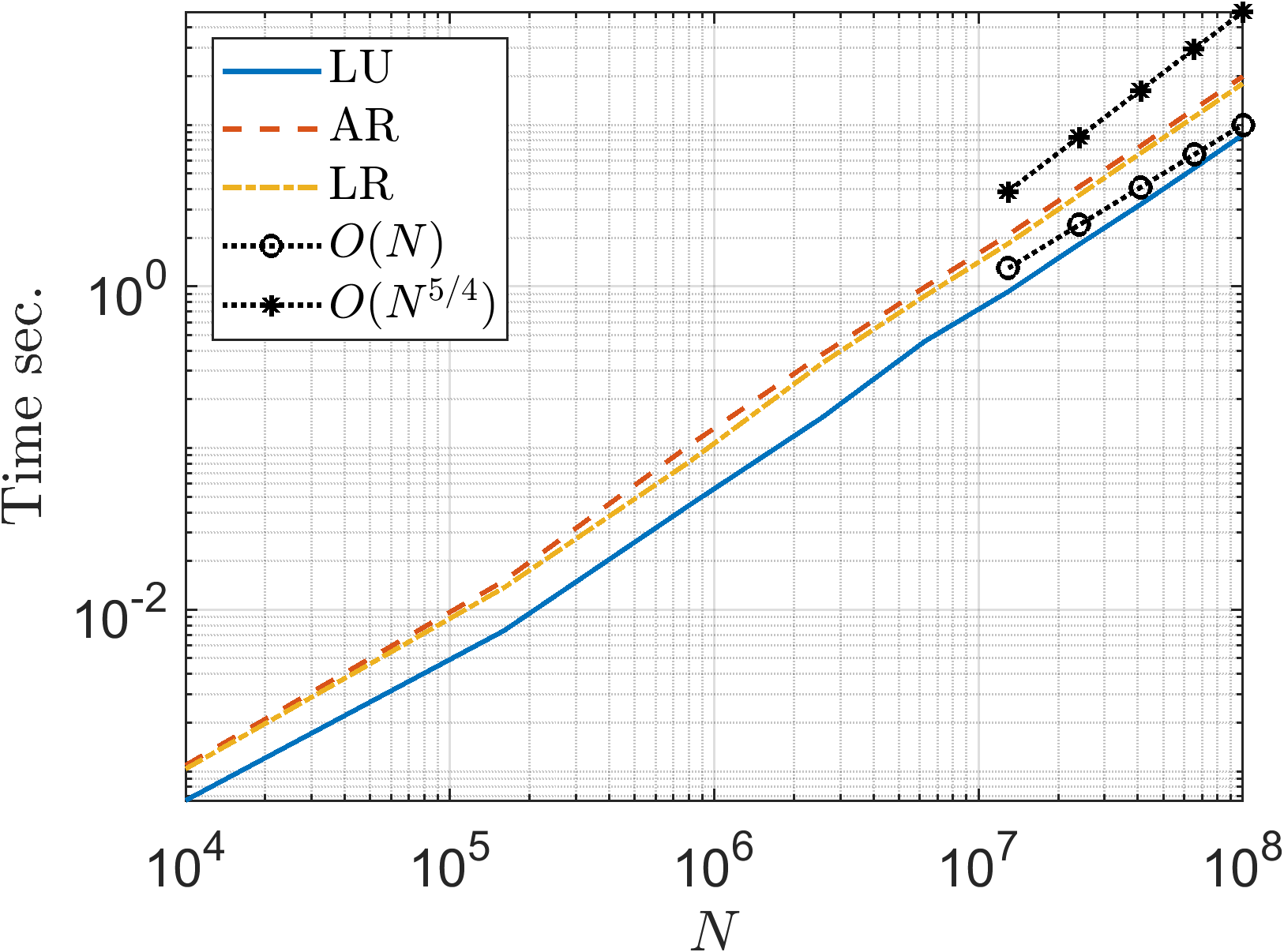}}  %\\
    \caption{Application costs for $N_s = N_t^3$.}
    \label{fig:application_plots}
\end{figure}

In Figure~\ref{fig:setup_plots} we show the setup time versus the number of degrees of freedom $N$ of the system, while in Figure~\ref{fig:application_plots} we show the application time of the three different methods studied. 
Notice that the time taken for setup is asymptotically $O(N)$ for all methods, but that the LU method shows an upward shift as the polynomial degree - and the bandwidth - increases. The total theoretical cost of the application for the three methods is $O(N^{5/4})$, as expected from the matrix vector multiplication in Step~3 and Step~5 of Algorithm~\ref{alg:gen_fast_diagonalization}.
However, numerical tests reveal a near-optimal reduced cost of order $O(N)$, and this is likely due to the high efficiency of the BLAS level 3 routines that perform the computational core of the application.
Finally Tables~\ref{tab:p1_total_times}--\ref{tab:p5_total_times} report the total amount of computational time required for the setup and application of the three methods. 
When the polynomial degree grows, the bandwidth of the matrices $\over{T}(\lambda)$ grows as well, thus the overall performance of the LU method worsens with respect to the AR or LR methods. 
We also see the performance of AR and LR is equivalent both in the setup and in the application cost.

\begin{table} 
	\centering
    \setlength{\tabcolsep}{4pt}  
    \begin{tabular}{c!{\vrule width 1pt}c|c|c|c|c|c!{\vrule width 1pt}c}
    \svhline              
    \hbox{\diagbox[linewidth=2\arrayrulewidth]{$N_s$}{$N_t$}  }
    & $16$ & $32$ & $64$ & $128$ & $256$ & $512$  & method \\     
    \svhline\noalign{\smallskip}                
                        & $ 3.7\cdot 10^{-3} $ & $ 7.4\cdot 10^{-3} $ & $ 1.3\cdot 10^{-2} $ & $ 2.6\cdot 10^{-2} $ & $ 4.7\cdot 10^{-2} $ & $ 9.5\cdot 10^{-2} $ & LU\\
    %\hline   
    $10^{3}$            & $ 3.0\cdot 10^{-3} $ & $ 5.2\cdot 10^{-3} $ & $ 1.1\cdot 10^{-2} $ & $ 3.9\cdot 10^{-2} $ & $ 1.1\cdot 10^{-1} $ & $ 5.3\cdot 10^{-1} $ & AR\\      
    %\hline   
                        & $ 3.2\cdot 10^{-3} $ & $ 5.3\cdot 10^{-3} $ & $ 1.0\cdot 10^{-2} $ & $ 3.8\cdot 10^{-2} $ & $ 1.0\cdot 10^{-1} $ & $ 5.2\cdot 10^{-1} $ & LR\\ 
    \noalign{\smallskip}\svhline\noalign{\smallskip}            
                        & $ 8.4\cdot 10^{-2} $ & $ 1.6\cdot 10^{-1} $ & $ 3.6\cdot 10^{-1} $ & $ 7.3\cdot 10^{-1} $ & $ 1.5\cdot 10^{0\phantom{-}}  $ & $ 3.0\cdot 10^{0\phantom{-}} $ & LU\\     
    %\hline   
    $2.7\cdot 10^4$     & $ 5.8\cdot 10^{-2} $ & $ 1.3\cdot 10^{-1} $ & $ 2.9\cdot 10^{-1} $ & $ 6.5\cdot 10^{-1} $ & $ 1.6\cdot 10^{0\phantom{-}}  $ & $ 4.3\cdot 10^{0\phantom{-}} $ & AR\\      
    %\hline   
                        & $ 6.5\cdot 10^{-2} $ & $ 1.3\cdot 10^{-1} $ & $ 2.6\cdot 10^{-1} $ & $ 6.3\cdot 10^{-1} $ & $ 1.6\cdot 10^{0\phantom{-}}  $ & $ 4.2\cdot 10^{0\phantom{-}} $ & LR\\ 
    \noalign{\smallskip}\svhline\noalign{\smallskip}            
                        & $ 4.2\cdot 10^{-1} $ & $ 8.7\cdot 10^{-1} $ & $ 1.8\cdot 10^{0\phantom{-}}  $ & $ 3.6\cdot 10^{0\phantom{-}}  $ & $ 7.1\cdot 10^{0\phantom{-}}  $ & $ 1.4\cdot 10^{1\phantom{-}} $ & LU\\     
    %\hline   
    $1.3\cdot 10^5$     & $ 3.4\cdot 10^{-1} $ & $ 7.0\cdot 10^{-1} $ & $ 1.5\cdot 10^{0\phantom{-}}  $ & $ 3.2\cdot 10^{0\phantom{-}}  $ & $ 7.4\cdot 10^{0\phantom{-}}  $ & $ 1.9\cdot 10^{1\phantom{-}} $ & AR\\      
    %\hline   
                        & $ 3.6\cdot 10^{-1} $ & $ 7.4\cdot 10^{-1} $ & $ 1.5\cdot 10^{0\phantom{-}}  $ & $ 3.1\cdot 10^{0\phantom{-}}  $ & $ 7.3\cdot 10^{0\phantom{-}}  $ & $ 1.9\cdot 10^{1\phantom{-}} $ & LR\\ 
    \noalign{\smallskip}\svhline\noalign{\smallskip}         
\end{tabular}
    \caption{Total computational time spent for solving the linear system with $p = 1$.}
    \label{tab:p1_total_times}
\end{table}

\begin{table} 
	\centering
    \setlength{\tabcolsep}{4pt}  
    \begin{tabular}{c!{\vrule width 1pt}c|c|c|c|c|c!{\vrule width 1pt}c}
    \svhline              
    \hbox{\diagbox[linewidth=2\arrayrulewidth]{$N_s$}{$N_t$}  }                
    & $16$ & $32$ & $64$ & $128$ & $256$ & $512$  & method \\     
    \svhline\noalign{\smallskip}                
                        & $ 4.8\cdot 10^{-3} $ & $ 7.8\cdot 10^{-3} $ & $ 1.6\cdot 10^{-2} $ & $ 3.3\cdot 10^{-2} $ & $ 6.7\cdot 10^{-2} $ & $ 1.4\cdot 10^{-1} $ & LU\\
    %\hline   
    $10^{3}$            & $ 4.5\cdot 10^{-3} $ & $ 4.7\cdot 10^{-3} $ & $ 1.0\cdot 10^{-2} $ & $ 2.9\cdot 10^{-2} $ & $ 1.0\cdot 10^{-1} $ & $ 5.0\cdot 10^{-1} $ & AR\\      
    %\hline   
                        & $ 5.1\cdot 10^{-3} $ & $ 4.8\cdot 10^{-3} $ & $ 1.0\cdot 10^{-2} $ & $ 2.7\cdot 10^{-2} $ & $ 1.0\cdot 10^{-1} $ & $ 5.0\cdot 10^{-1} $ & LR\\ 
    \noalign{\smallskip}\svhline\noalign{\smallskip}            
                        & $ 1.2\cdot 10^{-1} $ & $ 2.6\cdot 10^{-1} $ & $ 5.3\cdot 10^{-1} $ & $ 1.1\cdot 10^{0\phantom{-}} $ & $ 2.2\cdot 10^{0\phantom{-}}  $ & $ 4.4\cdot 10^{0\phantom{-}} $ & LU\\     
    %\hline   
    $2.7\cdot 10^4$     & $ 5.7\cdot 10^{-2} $ & $ 1.3\cdot 10^{-1} $ & $ 2.8\cdot 10^{-1} $ & $ 6.5\cdot 10^{-1          } $ & $ 1.6\cdot 10^{0\phantom{-}}  $ & $ 4.2\cdot 10^{0\phantom{-}} $ & AR\\      
    %\hline   
                        & $ 6.5\cdot 10^{-2} $ & $ 1.3\cdot 10^{-1} $ & $ 2.6\cdot 10^{-1} $ & $ 6.3\cdot 10^{-1          } $ & $ 1.5\cdot 10^{0\phantom{-}}  $ & $ 4.1\cdot 10^{0\phantom{-}} $ & LR\\ 
    \noalign{\smallskip}\svhline\noalign{\smallskip}            
                        & $ 6.1\cdot 10^{-1} $ & $ 1.2\cdot 10^{0\phantom{-}} $ & $ 2.6\cdot 10^{0\phantom{-}}  $ & $ 5.1\cdot 10^{0\phantom{-}}  $ & $ 1.0\cdot 10^{1\phantom{-}}  $ & $ 2.1\cdot 10^{1\phantom{-}} $ & LU\\     
    %\hline   
    $1.3\cdot 10^5$     & $ 3.4\cdot 10^{-1} $ & $ 7.0\cdot 10^{-1          } $ & $ 1.5\cdot 10^{0\phantom{-}}  $ & $ 3.2\cdot 10^{0\phantom{-}}  $ & $ 7.4\cdot 10^{0\phantom{-}}  $ & $ 1.9\cdot 10^{1\phantom{-}} $ & AR\\      
    %\hline   
                        & $ 3.6\cdot 10^{-1} $ & $ 7.2\cdot 10^{-1          } $ & $ 1.5\cdot 10^{0\phantom{-}}  $ & $ 3.1\cdot 10^{0\phantom{-}}  $ & $ 7.3\cdot 10^{0\phantom{-}}  $ & $ 1.9\cdot 10^{1\phantom{-}} $ & LR\\ 
    \noalign{\smallskip}\svhline\noalign{\smallskip}            
\end{tabular}
    \caption{Total computational time spent for solving the linear system with $p = 2$.}                          
    \label{tab:p2_total_times}            
\end{table}

\begin{table} 
	\centering
    \setlength{\tabcolsep}{4pt}  
    \begin{tabular}{c!{\vrule width 1pt}c|c|c|c|c|c!{\vrule width 1pt}c}
    \svhline              
    \hbox{\diagbox[linewidth=2\arrayrulewidth]{$N_s$}{$N_t$}  }
        & $16$ & $32$ & $64$ & $128$ & $256$ & $512$  & method \\     
    \svhline\noalign{\smallskip}                
                        & $ 5.2\cdot 10^{-3} $ & $ 1.0\cdot 10^{-3} $ & $ 2.1\cdot 10^{-2} $ & $ 4.4\cdot 10^{-2} $ & $ 9.5\cdot 10^{-2} $ & $ 1.9\cdot 10^{-1} $ & LU\\
    %\hline   
    $10^{3}$            & $ 2.6\cdot 10^{-3} $ & $ 6.6\cdot 10^{-3} $ & $ 1.0\cdot 10^{-2} $ & $ 2.9\cdot 10^{-2} $ & $ 1.0\cdot 10^{-1} $ & $ 4.8\cdot 10^{-1} $ & AR\\      
    %\hline   
                        & $ 2.9\cdot 10^{-3} $ & $ 8.2\cdot 10^{-3} $ & $ 1.0\cdot 10^{-2} $ & $ 2.7\cdot 10^{-2} $ & $ 9.9\cdot 10^{-2} $ & $ 4.7\cdot 10^{-1} $ & LR\\ 
    \noalign{\smallskip}\svhline\noalign{\smallskip}            
                        & $ 1.7\cdot 10^{-1} $ & $ 3.3\cdot 10^{-1} $ & $ 6.9\cdot 10^{-1} $ & $ 1.4\cdot 10^{0\phantom{-}} $ & $ 3.0\cdot 10^{0\phantom{-}}  $ & $ 5.9\cdot 10^{0\phantom{-}} $ & LU\\     
    %\hline   
    $2.7\cdot 10^4$     & $ 6.0\cdot 10^{-2} $ & $ 1.3\cdot 10^{-1} $ & $ 2.8\cdot 10^{-1} $ & $ 6.5\cdot 10^{-1          } $ & $ 1.6\cdot 10^{0\phantom{-}}  $ & $ 4.2\cdot 10^{0\phantom{-}} $ & AR\\      
    %\hline   
                        & $ 6.8\cdot 10^{-2} $ & $ 1.3\cdot 10^{-1} $ & $ 2.6\cdot 10^{-1} $ & $ 6.3\cdot 10^{-1          } $ & $ 1.5\cdot 10^{0\phantom{-}}  $ & $ 4.1\cdot 10^{0\phantom{-}} $ & LR\\ 
    \noalign{\smallskip}\svhline\noalign{\smallskip}            
                        & $ 7.9\cdot 10^{-1} $ & $ 1.6\cdot 10^{0\phantom{-}} $ & $ 3.5\cdot 10^{0\phantom{-}}  $ & $ 6.9\cdot 10^{0\phantom{-}}  $ & $ 1.4\cdot 10^{1\phantom{-}}  $ & $ 2.8\cdot 10^{1\phantom{-}} $ & LU\\     
    %\hline   
    $1.3\cdot 10^5$     & $ 3.4\cdot 10^{-1} $ & $ 7.0\cdot 10^{-1          } $ & $ 1.5\cdot 10^{0\phantom{-}}  $ & $ 3.2\cdot 10^{0\phantom{-}}  $ & $ 7.4\cdot 10^{0\phantom{-}}  $ & $ 1.9\cdot 10^{1\phantom{-}} $ & AR\\      
    %\hline   
                        & $ 3.9\cdot 10^{-1} $ & $ 7.3\cdot 10^{-1          } $ & $ 1.5\cdot 10^{0\phantom{-}}  $ & $ 3.1\cdot 10^{0\phantom{-}}  $ & $ 7.3\cdot 10^{0\phantom{-}}  $ & $ 1.9\cdot 10^{1\phantom{-}} $ & LR\\ 
    \noalign{\smallskip}\svhline\noalign{\smallskip}            
\end{tabular}
    \caption{Total computational time spent for solving the linear system with $p = 3$.}                          
    \label{tab:p3_total_times}            
\end{table}

\begin{table} 
	\centering
    \setlength{\tabcolsep}{4pt}  
    \begin{tabular}{c!{\vrule width 1pt}c|c|c|c|c|c!{\vrule width 1pt}c}
    \svhline              
    \hbox{\diagbox[linewidth=2\arrayrulewidth]{$N_s$}{$N_t$}  }
    & $16$ & $32$ & $64$ & $128$ & $256$ & $512$  & method \\     
    \svhline\noalign{\smallskip}                
                        & $ 6.2\cdot 10^{-3} $ & $ 1.3\cdot 10^{-2} $ & $ 2.7\cdot 10^{-2} $ & $ 5.6\cdot 10^{-2} $ & $ 1.2\cdot 10^{-1} $ & $ 2.5\cdot 10^{-1} $ & LU\\
    %\hline   
    $10^{3}$            & $ 2.6\cdot 10^{-3} $ & $ 4.7\cdot 10^{-3} $ & $ 1.0\cdot 10^{-2} $ & $ 2.9\cdot 10^{-2} $ & $ 1.0\cdot 10^{-1} $ & $ 4.8\cdot 10^{-1} $ & AR\\      
    %\hline   
                        & $ 2.8\cdot 10^{-3} $ & $ 4.8\cdot 10^{-3} $ & $ 1.0\cdot 10^{-2} $ & $ 2.7\cdot 10^{-2} $ & $ 9.7\cdot 10^{-2} $ & $ 4.7\cdot 10^{-1} $ & LR\\ 
    \noalign{\smallskip}\svhline\noalign{\smallskip}            
                        & $ 2.0\cdot 10^{-1} $ & $ 4.2\cdot 10^{-1} $ & $ 8.8\cdot 10^{-1} $ & $ 1.8\cdot 10^{0\phantom{-}} $ & $ 3.7\cdot 10^{0\phantom{-}}  $ & $ 7.4\cdot 10^{0\phantom{-}} $ & LU\\     
    %\hline   
    $2.7\cdot 10^4$     & $ 5.9\cdot 10^{-2} $ & $ 1.3\cdot 10^{-1} $ & $ 2.8\cdot 10^{-1} $ & $ 6.5\cdot 10^{-1          } $ & $ 1.6\cdot 10^{0\phantom{-}}  $ & $ 4.2\cdot 10^{0\phantom{-}} $ & AR\\      
    %\hline   
                        & $ 6.5\cdot 10^{-2} $ & $ 1.3\cdot 10^{-1} $ & $ 2.6\cdot 10^{-1} $ & $ 6.3\cdot 10^{-1          } $ & $ 1.5\cdot 10^{0\phantom{-}}  $ & $ 4.1\cdot 10^{0\phantom{-}} $ & LR\\ 
    \noalign{\smallskip}\svhline\noalign{\smallskip}            
                        & $ 9.6\cdot 10^{-1} $ & $ 2.1\cdot 10^{0\phantom{-}} $ & $ 4.3\cdot 10^{0\phantom{-}}  $ & $ 8.6\cdot 10^{0\phantom{-}}  $ & $ 1.8\cdot 10^{1\phantom{-}}  $ & $ 3.5\cdot 10^{1\phantom{-}} $ & LU\\     
    %\hline   
    $1.3\cdot 10^5$     & $ 3.4\cdot 10^{-1} $ & $ 6.9\cdot 10^{-1          } $ & $ 1.5\cdot 10^{0\phantom{-}}  $ & $ 3.2\cdot 10^{0\phantom{-}}  $ & $ 7.4\cdot 10^{0\phantom{-}}  $ & $ 1.9\cdot 10^{1\phantom{-}} $ & AR\\      
    %\hline   
                        & $ 3.5\cdot 10^{-1} $ & $ 7.1\cdot 10^{-1          } $ & $ 1.4\cdot 10^{0\phantom{-}}  $ & $ 3.1\cdot 10^{0\phantom{-}}  $ & $ 7.2\cdot 10^{0\phantom{-}}  $ & $ 1.8\cdot 10^{1\phantom{-}} $ & LR\\ 
    \noalign{\smallskip}\svhline\noalign{\smallskip}            
\end{tabular}
    \caption{Total computational time spent for solving the linear system with $p = 4$.}                          
    \label{tab:p4_total_times}            
\end{table}

\begin{table} 
	\centering
    \setlength{\tabcolsep}{4pt}  
    \begin{tabular}{c!{\vrule width 1pt}c|c|c|c|c|c!{\vrule width 1pt}c}
    \svhline              
    \hbox{\diagbox[linewidth=2\arrayrulewidth]{$N_s$}{$N_t$}  }
    & $16$ & $32$ & $64$ & $128$ & $256$ & $512$  & method \\     
    \svhline\noalign{\smallskip}        
                        & $ 8.3\cdot 10^{-3} $ & $ 1.7\cdot 10^{-2} $ & $ 3.3\cdot 10^{-2} $ & $ 7.2\cdot 10^{-2} $ & $ 1.5\cdot 10^{-1} $ & $ 3.1\cdot 10^{-1} $ & LU\\
    %\hline   
    $10^{3}$            & $ 2.7\cdot 10^{-3} $ & $ 4.7\cdot 10^{-3} $ & $ 1.0\cdot 10^{-2} $ & $ 2.9\cdot 10^{-2} $ & $ 1.0\cdot 10^{-1} $ & $ 4.8\cdot 10^{-1} $ & AR\\      
    %\hline   
                        & $ 2.9\cdot 10^{-3} $ & $ 4.7\cdot 10^{-3} $ & $ 1.1\cdot 10^{-2} $ & $ 2.7\cdot 10^{-2} $ & $ 9.5\cdot 10^{-2} $ & $ 4.7\cdot 10^{-1} $ & LR\\ 
    \noalign{\smallskip}\svhline\noalign{\smallskip}            
                        & $ 2.4\cdot 10^{-1} $ & $ 5.1\cdot 10^{-1} $ & $ 1.1\cdot 10^{0\phantom{-}} $ & $ 2.2\cdot 10^{0\phantom{-}} $ & $ 4.6\cdot 10^{0\phantom{-}}  $ & $ 9.1\cdot 10^{0\phantom{-}} $ & LU\\     
    %\hline   
    $2.7\cdot 10^4$     & $ 5.6\cdot 10^{-2} $ & $ 1.3\cdot 10^{-1} $ & $ 2.8\cdot 10^{-1          } $ & $ 6.4\cdot 10^{-1          } $ & $ 1.5\cdot 10^{0\phantom{-}}  $ & $ 4.2\cdot 10^{0\phantom{-}} $ & AR\\      
    %\hline   
                        & $ 6.2\cdot 10^{-2} $ & $ 1.3\cdot 10^{-1} $ & $ 2.6\cdot 10^{-1          } $ & $ 6.1\cdot 10^{-1          } $ & $ 1.5\cdot 10^{0\phantom{-}}  $ & $ 4.1\cdot 10^{0\phantom{-}} $ & LR\\ 
    \noalign{\smallskip}\svhline\noalign{\smallskip}            
                        & $ 1.2\cdot 10^{0\phantom{-}} $ & $ 2.5\cdot 10^{0\phantom{-}} $ & $ 5.3\cdot 10^{0\phantom{-}}  $ & $ 1.1\cdot 10^{1\phantom{-}}  $ & $ 2.1\cdot 10^{1\phantom{-}}  $ & $ 4.3\cdot 10^{1\phantom{-}} $ & LU\\     
    %\hline   
    $1.3\cdot 10^5$     & $ 3.4\cdot 10^{-1          } $ & $ 6.8\cdot 10^{-1          } $ & $ 1.4\cdot 10^{0\phantom{-}}  $ & $ 3.1\cdot 10^{0\phantom{-}}  $ & $ 7.3\cdot 10^{0\phantom{-}}  $ & $ 1.9\cdot 10^{1\phantom{-}} $ & AR\\      
    %\hline   
                        & $ 3.5\cdot 10^{-1          } $ & $ 7.1\cdot 10^{-1          } $ & $ 1.4\cdot 10^{0\phantom{-}}  $ & $ 3.0\cdot 10^{0\phantom{-}}  $ & $ 7.1\cdot 10^{0\phantom{-}}  $ & $ 1.8\cdot 10^{1\phantom{-}} $ & LR\\ 
    \noalign{\smallskip}\svhline\noalign{\smallskip}            
\end{tabular}
    \caption{Total computational time spent for solving the linear system with $p = 5$.}                          
    \label{tab:p5_total_times}            
\end{table}

\subsection{Application as a preconditioner}\label{sec:numeric-preconditioner}
\edef\psd{\widehat{\Omega}_{s}}
In this section, we apply the LU, AR, and LR methods as  left preconditioners in the context of isogeometric discretizations. 
In isogeometric analysis (IGA), introduced in \cite{MR2152382},
the space domain $\sd$ is the
image through a bi-Lipschitz map $G$ of a cuboid,
usually $\psd:=[0,1]^d$.
The discrete space defined on $\sd\times \td$ is $V:=V_s\otimes V_t$ where $V_s$ is the push-forward 
of a tensor product spline space
$\widehat V_s$ defined on $\psd$
\[
V_s=\{\hat v(G^{-1}(\mbf x)): \hat v\in \widehat V_s\},
\]
see \cite{MR3202239} for more details, and $V_t$ is a spline space on $\td$.

The domain $\sd$, in general, is not a Cartesian product
and consequently the stiffness matrix $\mbf A_s$ and the Gram
matrix $\mbf M_s$ for the space $V_s$ do not satisfy \eqref{eq:Ms-tp}
nor \eqref{eq:kronecker_structure_of_stiffness_matrix}
and the generalized eigendecomposition of the pair 
$(\mbf A_s,\mbf M_s)$ is computationally infeasible.

Nevertheless the stiffness matrix $\widehat{\mbf A}_s$ and the Gram
matrix $\widehat{\mbf M}_s$ constructed for the space $\widehat V_s$   on the Cartesian domain $\widehat{\Omega}_s$
satisfy both conditions so that the presented methods (LU, AR, LR) can be
efficiently applied to  the heat operator on $\psd\times\td$ that is
\[
\widehat{\mbf A} = \mbf{A}_t \otimes \widehat{\mbf{M}}_s + \mbf{M}_t \otimes \widehat{\mbf A}_s.
\]
 
The matrices $\mbf M_s$ and $\widehat{\mbf M}_s$ are spectrally
equivalent with constants that depend only on $G$ \cite{loli2022easy}, similarly 
$\mbf A_s$ and $\widehat{\mbf A}_s$ \cite{igaFD2016}.
As in \cite{loli2020efficient} we use $\widehat{\mbf A}$ as
a   left preconditioner , i.e. we solve the system
\[
\widehat{\mbf A}^{-1} \mbf A \mbf u =\widehat{\mbf A}^{-1} \mbf f,
\]
where the application of $\widehat{\mbf A}^{-1}$ to a vector is performed by the described LU, AR and LR methods.
% to accelerate an iterative Krylov solver.

\begin{figure}[t]
    \centering
    \subfloat[][Cube.\label{fig:cube}]
    {\includegraphics[trim=5cm 0cm 5cm 0cm, clip=true, scale=0.25]{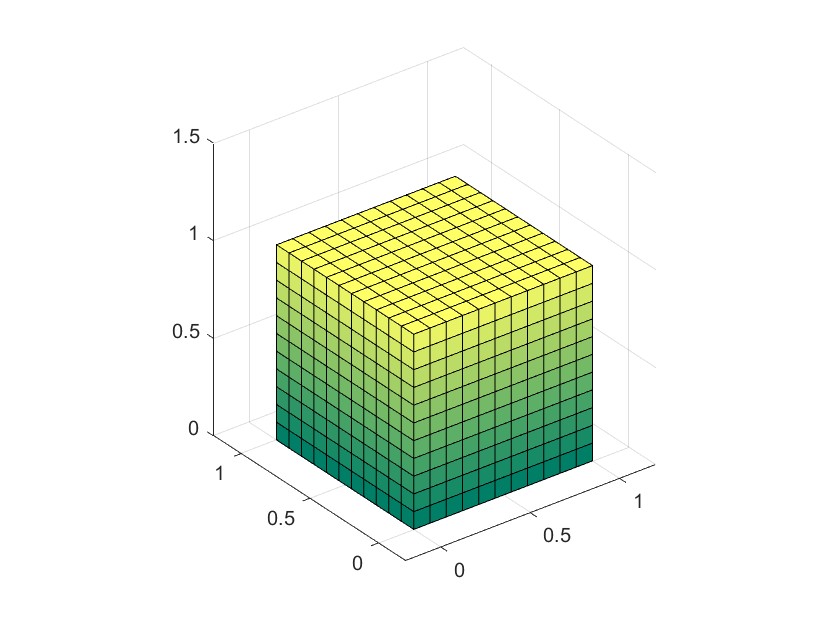}}  \quad 
    \subfloat[][Rotated quarter of annulus.\label{fig:rev-quarter}]
    {\includegraphics[trim=5cm 0cm 5cm 0cm, clip=true, scale=0.25]{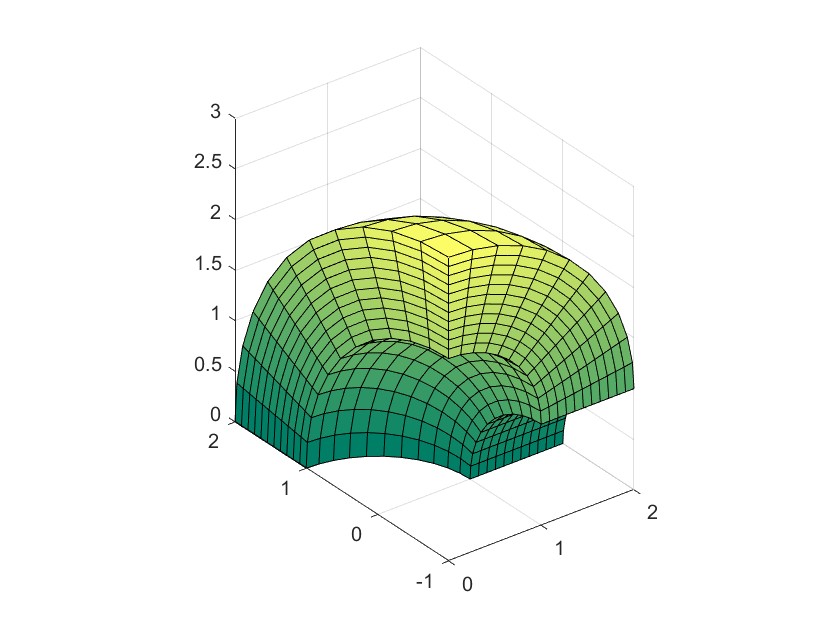}}  %\\
    \caption{Computational domains $\sd$.}
    \label{fig:isogeometric-domain}
\end{figure}

In the presented tests $G$ maps the cube in Figure~\ref{fig:cube}
to the geometry in Figure \ref{fig:rev-quarter}
that is a rotational solid obtained by rotating by $\frac{\pi}{2}$ 
a 2D quarter of annulus, with center in the origin
internal radius 1 and external radius 2,
around the axis $\{(x_1,-1,0) | x_1 \in \R\}$.
The final time is $T=1$, Dirichlet and initial conditions are such that 
$u(\vect{x},t) = -(x_1^2 + x_2^2-1)(x_1^2+x_2^2-4)x_1x_2^2 \sin(x_3) \sin(t)$ 
is the exact solution. %with constant coefficients $\gamma = \nu = 1$.
The discrete space $\widehat V$ contains splines of the same degree in space and in time $p_t=p_s=:p$, uniform knot vectors and $N_{s,l}=N_t$, $l=1,2,3$.
\begin{table}[h]
    \centering
    \setlength{\tabcolsep}{4pt}
    \begin{tabular}{c!{\vrule width 1pt}c|c|c|c|c}
    \svhline              
    \hbox{\diagbox[linewidth=2\arrayrulewidth]{$N_t$}{$p$}} 
    & $1$ &$2$ & $3$ & $4$ & $5$   \\
    \svhline\noalign{\smallskip}
    8  & $37 $  &  $38 $  &  $41 $  &   $43$ &   $45 $ \\
    16 & $46 $  &  $49 $  &  $52 $  &   $55$ &   $58 $ \\
    32 & $53 $  &  $56 $  &  $60 $  &   $63$ &   $65 $ \\
    64 & $59 $  &  $62 $  &  $64 $  &   $67$ &   $70 $ \\
    \noalign{\smallskip}\svhline\noalign{\smallskip}
    \end{tabular}
    \caption{Number of iterations in  GMRES solver with preconditioner $\widehat{\mbf A}$}.
    \label{tab:preconditioner-performance-iterations}
\end{table}

\begin{table}[h]
    \centering
    \setlength{\tabcolsep}{4pt}
    \begin{tabular}{c!{\vrule width 1pt}c|c|c|c|c!{\vrule width 1pt}c}
    \svhline              
    \hbox{\diagbox[linewidth=2\arrayrulewidth]{$N_t$}{$p$}} 
    & $1$ &$2$ & $3$ & $4$ & $5$ & method  \\
    \svhline\noalign{\smallskip}
       & $ 1.7\cdot 10^{-1          } $  &  $ 1.1\cdot 10^{-1          } $  &  $ 9.0\cdot 10^{-2          } $  &   $ 1.4\cdot 10^{-1          } $  &   $ 1.1\cdot 10^{0} $ & LU\\
    8  & $ 2.3\cdot 10^{-1          } $  &  $ 1.3\cdot 10^{-1          } $  &  $ 1.9\cdot 10^{-1          } $  &   $ 2.4\cdot 10^{-1          } $  &   $ 3.2\cdot 10^{0} $ & AR\\
       & $ 1.4\cdot 10^{-1          } $  &  $ 1.2\cdot 10^{-1          } $  &  $ 1.8\cdot 10^{-1          } $  &   $ 2.4\cdot 10^{-1          } $  &   $ 3.2\cdot 10^{0} $ & LR\\
       \noalign{\smallskip} \svhline \noalign{\smallskip} 
       & $ 6.1\cdot 10^{-1          } $  &  $ 1.0\cdot 10^{0\phantom{-}} $  &  $ 2.0\cdot 10^{0\phantom{-}} $  &   $ 3.6\cdot 10^{0\phantom{-}} $  &   $ 6.2\cdot 10^{0} $ & LU\\
    16 & $ 1.9\cdot 10^{0\phantom{-}} $  &  $ 2.5\cdot 10^{0\phantom{-}} $  &  $ 4.1\cdot 10^{0\phantom{-}} $  &   $ 6.6\cdot 10^{0\phantom{-}} $  &   $ 9.9\cdot 10^{0} $ & AR\\
       & $ 1.8\cdot 10^{0\phantom{-}} $  &  $ 2.5\cdot 10^{0\phantom{-}} $  &  $ 4.1\cdot 10^{0\phantom{-}} $  &   $ 6.5\cdot 10^{0\phantom{-}} $  &   $ 9.9\cdot 10^{0} $ & LR\\
       \noalign{\smallskip} \svhline \noalign{\smallskip} 
       & $ 1.4\cdot 10^{1\phantom{-}} $  &  $ 2.3\cdot 10^{1\phantom{-}} $  &  $ 4.6\cdot 10^{1\phantom{-}} $  &   $ 8.3\cdot 10^{1\phantom{-}} $  &   $ 1.6\cdot 10^{2} $ & LU\\
    32 & $ 4.2\cdot 10^{1\phantom{-}} $  &  $ 5.9\cdot 10^{1\phantom{-}} $  &  $ 9.5\cdot 10^{1\phantom{-}} $  &   $ 1.6\cdot 10^{2\phantom{-}} $  &   $ 2.7\cdot 10^{2} $ & AR\\
       & $ 4.2\cdot 10^{1\phantom{-}} $  &  $ 5.8\cdot 10^{1\phantom{-}} $  &  $ 9.3\cdot 10^{1\phantom{-}} $  &   $ 1.7\cdot 10^{2\phantom{-}} $  &   $ 2.7\cdot 10^{2} $ & LR\\
       \noalign{\smallskip} \svhline \noalign{\smallskip} 
       & $ 3.3\cdot 10^{2\phantom{-}} $  &  $ 5.1\cdot 10^{2\phantom{-}} $  &  $ 8.6\cdot 10^{2\phantom{-}} $  &   $ 1.8\cdot 10^{3\phantom{-}} $  &   $ 3.1\cdot 10^{3} $ & LU\\
    64 & $ 1.1\cdot 10^{3\phantom{-}} $  &  $ 1.4\cdot 10^{3\phantom{-}} $  &  $ 2.0\cdot 10^{3\phantom{-}} $  &   $ 3.6\cdot 10^{3\phantom{-}} $  &   $ 5.7\cdot 10^{3} $ & AR\\
       & $ 1.1\cdot 10^{3\phantom{-}} $  &  $ 1.3\cdot 10^{3\phantom{-}} $  &  $ 2.0\cdot 10^{3\phantom{-}} $  &   $ 3.6\cdot 10^{3\phantom{-}} $  &   $ 5.5\cdot 10^{3} $ & LR\\
    \noalign{\smallskip}\svhline\noalign{\smallskip}
    \end{tabular}
    \caption{Time  used by the preconditioned GMRES solver, depending on the algorithm (LU, AR, LR), the degree $p$ and the number of elements in the time direction $N_t$}.
    \label{tab:preconditioner-performance-times}
\end{table}

The   preconditioned linear system is solved by GMRES, with tolerance equal to $10^{-8}$ and with the null vector as initial guess. 
Table~\ref{tab:preconditioner-performance-iterations} reports the number of iterations required to converge using $\widehat{\mbf A}$ as preconditioner for different choices of the degree $p$ and number of elements $N_t$.
We emphasise that all methods (LU, AR and LR) required the same number of iterations. This is not not surprising, since these methods are mathematically equivalent.
Finally, in Table \ref{tab:preconditioner-performance-times} we report for the three methods the total solving time, % as presented in \cite[Table 4]{loli2020efficient}, 
showing that LU method allows for a faster computation of the solution with respect to AR and LR methods that are equivalent.

\section{Conclusions}\label{sec:conclusions}
%quale era il problema...
In this work we have investigated several Kronecker-based methods for solving linear evolution problems 
with particular emphasis on the heat equation. 
The Diagonalization in Time (DT) %an extension of the highly efficient fast diagonalization (FD) method for elliptic problems,
 exhibits significant instabilities when applied to the heat equation.
%cosa abbiamo fatto
To address these limitations, we examined three alternative approaches: LU~factorization (LU), arrowhead factorization (AR), and low-rank modification (LR). 
These methods replace the unstable generalized eigendecomposition in time \eqref{eq:eigen-decomposition-t} with more robust factorizations, 
thereby maintaining computational efficiency while ensuring stability in the solution of the heat equation.
These methods require a tensor product structure of the discrete space
when used as direct solvers. When this assumption fails, but it is satisfied by
an approximate linear system, they can be used as preconditioners as demonstrated in one of the presented numerical experiments. 

\begin{table}
	\centering
    \begin{tabular}{c|cc|cc}
    \svhline\noalign{\smallskip}
    method & Step~1 & Step~2  & Step~3,5 & Step~4  
\\
    \noalign{\smallskip}\svhline\noalign{\smallskip} 
    DT 
    & $\sum_{l=1}^d N_{s,l}^3$ 
    & $N_t^3+\phantom{\rank^2}N$
    & $N (N_t + \sum_{l=1}^d N_{s,l})$
    & $\phantom{b}N$
\\  LU
    & $\sum_{l=1}^d N_{s,l}^3$
    & $\phantom{N_t^3+ \mbox{}}\band^2 N$
    & $ N (\ \phantom{N_t+} \sum_{l=1}^d N_{s,l} )$
    & $\band N$
\\  AR    
    & $\sum_{l=1}^d N_{s,l}^3$
    & $N_t^3+\phantom{\rank^2}N$
    & $ N ( N_t+\sum_{l=1}^d N_{s,l})$
    & $\phantom{b}N$
\\  LR
    & $\sum_{l=1}^d N_{s,l}^3$
    & $N_t^3+\rank^2 N$
    & $N( N_t+\sum_{l=1}^d N_{s,l})$
    & $\rank N$
\\
   \noalign{\smallskip}\svhline\noalign{\smallskip}            
    \end{tabular}
    \caption{The costs of setup and application for all the considered methods. Recall that $\band$ is the bandwidth of $\mbf{A}_t$ and $\mbf{M}_t$ and $\rank$ is the rank of the low-rank modification.}
    \label{tab:full-cost-table}
    \setlength{\tabcolsep}{12pt}
\end{table}
Recalling the theoretical asymptotic complexity of the different steps of the presented methods in Table~\ref{tab:full-cost-table} % (see Sections~\ref{sub_sec:banded}~-~\ref{sub_sub_sec:low_rank} for details).
we see: stable methods have the same theoretical cost with respect to DT method, Step~3 and Step~5 dominate the application phase and LU method has the fastest application since $\band\le N_t$. 
These findings are confirmed by the numerical tests, additionally showing that LU method has the slowest setup.

As the AR and LR methods can be parallelized in both spatial and temporal dimensions they are suitable for complex and large-scale simulations.
The LU method does not require that the pair $(\mbf A_t,\mbf M_t)$ has a generalized eigendecomposition and as such is applicable to a broader
class of operators acting on functions of time, allowing for generalization to a broader range of partial differential equations, e.g. the least square formulation of the Schrödinger equation in \cite{bressan2023space}.
Moreover, numerical tests in dimension $d=3$ reveal that the computational complexity of the methods studied is of order $O(N)$, 
despite the fact that the theoretically dominant computational cost is given by Step~3 and Step~5 in Table \ref{tab:full-cost-table}, that for $N_t = N_{s,l} = N ^{1/(d+1)}$ is $O(N^{(d+2)/(d+1)})$. 

%ricerca futura
Future research could build on these results by investigating their applicability to variable coefficient problems, cases beyond the parabolic framework, and nonlinear equations.
In conclusion, while these methods exhibit optimal behavior in academic tests, several challenges remain when addressing real-world applications.

\begin{acknowledgement}
The authors are members of the Gruppo Nazionale Calcolo Scientifico - Istituto Nazionale di Alta Matematica (GNCS-INdAM).
A. Kushova, G. Loli  and M. Montardini are partially supported by GNCS-INdAM project ‘‘Metodi numerici efficienti per la risoluzione di PDEs’’.
M. Tani  is partially supported by GNCS-INdAM project “Tecnologie spline efficienti per l'analisi isogeometrica e l'approssimazione di dati”.
A. Bressan is partially supported by GNCS-INdAM project "Metodi numerici avanzati per la poromeccanica: proprieta' teoriche ed aspetti computazionali”.
A. Kushova, G. Loli,  G. Sangalli and M. Tani acknowledge support from the MUR
through the PRIN 2022 PNRR project NOTES (No. P2022NC97R).
A. Bressan and M. Tani acknowledge support from the
MUR through the PRIN project COSMIC (No. 2022A79M75).
M. Montardini acknowledges support from the MUR through the PRIN 2022 PNRR project HEXAGON, Italy (No. P20227CTY3).
G. Loli and G. Sangalli acknowledge support from PNRR-M4C2-I1.4-NC-HPC-Spoke6, Italy. 

\end{acknowledgement}

%\section*{Appendix}
%If needed.

\bibliographystyle{plain}
\bibliography{references}
\end{document}